\documentclass[11pt]{article}
\usepackage[T1]{fontenc}
\usepackage[latin9]{inputenc}
\usepackage[top=3cm, bottom=3cm, left=1.2cm, right=1.2cm]{geometry}
\usepackage{amsmath}
\usepackage{amsthm}
\usepackage{amssymb}
\usepackage[unicode=true,pdfusetitle,
bookmarks=true,bookmarksnumbered=true,bookmarksopen=true,bookmarksopenlevel=3,
breaklinks=false,pdfborder={0 0 1},backref=false,colorlinks=false]
{hyperref}

\usepackage[english]{babel}
\usepackage{verbatim,url}            
\usepackage{booktabs}       
\usepackage{amsfonts}       
\usepackage{nicefrac}       
\usepackage[mathscr]{eucal}
\usepackage[pdftex]{graphicx}
\usepackage[colorinlistoftodos,bordercolor=orange,backgroundcolor=orange!20,linecolor=orange,textsize=scriptsize]{todonotes}
\usepackage{subfig}
\usepackage{color, makecell, afterpage}%

\makeatletter
\theoremstyle{plain}
\newtheorem{thm}{\protect\theoremname}
\theoremstyle{plain}

\theoremstyle{plain}
\newtheorem{prop}[thm]{\protect\propositionname}
\theoremstyle{definition}
\newtheorem*{remark}{Remark}
\usepackage{multirow}

\usepackage{lmodern}
\newcommand{\1}{\mbox{1\hspace{-1mm}I}}
\numberwithin{equation}{section}

\usepackage{calc}
\makeatother

\usepackage[round]{natbib}
\usepackage{cleveref}
\usepackage[short]{optidef}

\newcommand{\bE}{\mathbb{E}}

\newcommand{\RR}{\mathbb{R}}

\newcommand{\R}{\RR} 

\DeclareMathOperator*{\Id}{Id}
\newcommand{\cH}{{\mathcal H}}
\newcommand{\cL}{{\mathcal L}}

\newcommand{\cJ}{{\mathcal J}}
\newcommand{\cF}{{\mathcal F}}
\newcommand{\cP}{{\mathcal P}}
\newcommand{\cW}{{\mathcal W}}

\newcommand{\Vb}{{\bar V}}
\newcommand{\Mb}{{\bar M}}
\newcommand{\Xb}{{\bar X}}
\newcommand{\Fb}{{\bar F}}
\newcommand{\Ub}{{\bar U}}
\newcommand{\Yb}{{\bar Y}}
\newcommand{\Jb}{{\bar J}}

\newcommand{\Xt}{{\widetilde X}}
\newcommand{\Vt}{{\widetilde V}}
\newcommand{\Ut}{{\widetilde U}}
\newcommand{\Jt}{{\widetilde J}}

\providecommand{\noridx}[2]{\lVert{#1}\rVert_{#2}}

\DeclareMathOperator*{\Div}{div}
\DeclareMathOperator*{\KL}{KL}

\addto\captionsamerican{\renewcommand{\lemmaname}{Lemma}}
\addto\captionsamerican{\renewcommand{\propositionname}{Proposition}}
\addto\captionsamerican{\renewcommand{\theoremname}{Theorem}}
\addto\captionsenglish{\renewcommand{\lemmaname}{Lemma}}
\addto\captionsenglish{\renewcommand{\propositionname}{Proposition}}
\addto\captionsenglish{\renewcommand{\theoremname}{Theorem}}
\providecommand{\lemmaname}{Lemma}
\providecommand{\propositionname}{Proposition}
\providecommand{\theoremname}{Theorem}
\begin{document}
\global\long\def\1{\mbox{1\hspace{-1mm}I}}%

\title{REPRODUCING KERNEL APPROACH TO LINEAR QUADRATIC MEAN FIELD CONTROL PROBLEMS\\
}
\author{Pierre-Cyril Aubin-Frankowski\footnote{INRIA-Département d'Informatique de l'École Normale Supérieure, PSL, Research University, Paris, France }\\
Alain Bensoussan\footnote{International Center for Decision and Risk Analysis, Jindai School of Management, University of Texas at Dallas, School of Data Science, City University Hong Kong}
\\
}
\maketitle

\begin{abstract}
	Mean-field control problems have received continuous interest over the last decade. Despite being more intricate than in classical optimal control, the linear-quadratic setting can still be tackled through Riccati equations. Remarkably, we demonstrate that another significant attribute extends to the mean-field case: the existence of an intrinsic reproducing kernel Hilbert space associated with the problem. Our findings reveal that this Hilbert space not only encompasses deterministic controlled push-forward mappings but can also represent of stochastic dynamics. Specifically, incorporating Brownian noise affects the deterministic kernel through a conditional expectation, to make the trajectories adapted. Introducing reproducing kernels allows us to rewrite the mean-field control problem as optimizing over a Hilbert space of trajectories rather than controls. This framework even accommodates nonlinear terminal costs, without resorting to adjoint processes or Pontryagin's maximum principle, further highlighting the versatility of the proposed methodology.
\end{abstract}



\section{INTRODUCTION}

\paragraph{\bf Context.} In this article, we are interested in solving linear-quadratic (LQ) mean-field control (MFC) problems of the form
\begin{mini}{\substack{(\alpha_T)_{t\in[0,T]}\in \cP(\R^n)^{[0,T]},\\ V(\cdot)\in L^2([0,T]\times \R^n; \R^d)}}{\Phi(\mu_T)+\int_{0}^T \int_{\R^n}L(t,x,V(t,x))d\mu_t(x)dt}{\label{eq:obj_FP}}{}%
	\addConstraint{\partial_t \mu_t+\Div\left(g\left(t,x,V(t,x), \int_{\R^n}x'd\mu_t(x')\right)\cdot\mu_t\right)=0}{,\, \forall t\in[0,T]}{}
	\addConstraint{\mu_0=\mu^0}{}{}
\end{mini}
where $\mu^0\in \cP(\R^n)$ is a given probability measure over $\R^n$, $L$ is a sum of convex quadratic functions in $x$ and in $v$, and $g$ is linear in its arguments other then $t$, all coefficients being assumed deterministic. This problem is of the mean-field type because integrals over $\mu_t$ appear in both the dynamics and in the objective. Notably, the continuity equation can also be described by a linear controlled McKean--Vlasov process over a state in $\R^n$. This allows for several equivalent definitions of problem \eqref{eq:obj_FP} as explained in \cite{Bensoussan2013mfg} or \cite[Section 6.2]{Carmona2018}. We will focus on yet another reformulation, first presented in \cite{Bensoussan2020controlHilbert}, based on controlled push-forward maps. We first consider deterministic dynamics and $\Phi$ quadratic in the mean and variance of $\mu_t$. We later move to nonlinear $\Phi$ and stochastic dynamics.

\paragraph{\bf Related work.} Mean-field control originated from diverse perspectives such as the limit of specific multi-agent interactions as done in \cite{Fornasier2014}, or, when starting from controlled Fokker--Planck equations, as a control counterpart to mean-field games as summarized in \cite{Bensoussan2013mfg}, or as the study of controlled McKean--Vlasov dynamics as thoroughly done by \cite{Carmona2012,Carmona2015}. Both of the latter two approaches are related in spirit to the study of the maximum principle by \cite{Andersson2010}. The linear-quadratic (LQ) framework, akin to its classical counterpart, emerges as a natural starting point. Early instances of LQ problems can be traced back to \cite[Section 5]{Huang2007} for the game setting or in \cite[Section 5]{Andersson2010} for Mean-Variance Portfolio Selection. For a survey of the controlled McKean--Vlasov literature, we refer to \cite{Yong2013lq}, and to \cite{Pham2016} for the more intricate case of random coefficients. While adapting traditional solving techniques specific to classical LQ optimal control, predominantly through Riccati equations, is not straightforward, it has been done in \cite{Yong2013lq,Sun2020}, \cite[Section 6.3]{Bensoussan2013mfg}, \cite[Section 6]{Pham2017}, and for infinite-dimensional states in \cite{Lu2020sto}. However eventually one wants to allow for some nonlinearities in $\Phi$, such as entropy-like terms, and thus involve some differential calculus which goes beyond Riccati equations.

Several very interesting directions have been pursued to be able to differentiate in the context of measures spaces and continuity equations. The foremost, and the one we are most apart from, is Otto's calculus introduced in \cite{Otto2001}, where the set of probability measures is equipped with a Riemannian-like structure based on the Wasserstein distance. However this perspective was mostly explored from a calculus of variations perspective as in the reference book by \cite{ambrosio2008gradient} and has only recently been considered for control settings by \cite{Gangbo2015} on the torus, by \cite{Bongini2017} to show that the optimality of some systems corresponded to Hamiltonian flows and by \cite{Bonnet2020,Bonnet2022} for more general differential inclusions on $\R^n$. Most importantly, in Otto's calculus, the linearity occurs on the tangent spaces to the Wasserstein space, which are $L^2$-spaces of vector fields. These, when concatenated, provide push-forward maps $X_V^t$ for the measures, mapping $\mu^0$ to $\mu_t$. A nice perk of linear dynamics is that we can then introduce the Hilbert space of square-integrable time-varying push-forward maps $L^2(0,T;L_{\mu^0}^2(\R^n,\R^n))$ as first done by one of the authors in \cite[Section 3]{Bensoussan2020controlHilbert}. This differs from a procedure proposed by P.\-L.\ Lions, which ``flattens'' the Wasserstein space by lifting it to the Hilbert space of random variables with values in $\R^n$. We refer to \cite[Section 6.1]{Cardaliaguet2012} and \cite[Section 3.1]{Carmona2015} for the consequent differentiation rules of this other choice of lifting.

In a distinct context of uncovering connections between optimal control and reproducing kernels, we showed recently in \cite{aubin2020hard_control} that, when the focus of the problem is on the controlled trajectories --here the role played by the push-forwards-- classical LQ optimal control problems had a previously unexplored reproducing kernel Hilbert space (RKHS) structure, allowing to tackle easily state constraints. In the follow-up study \cite{aubin2022operator}, we extended the approach to infinite-dimensional settings, i.e.\ control of parabolic equations such as the heat equation. We here show that the kernel viewpoint also applies to LQ MFC. It is worth noting that using matrix-valued kernels in the context of push-forward maps has been previously explored by \cite{Arguillre2015} for shape deformations. However \cite{Arguillre2015} used off-the-shelf kernels, like the Gaussian kernel, defined over the space variable and encoding the control vector fields $V(\cdot)$ of the dynamics ($\dot{x}=v(s)$). Our operator-valued kernels are defined instead on the time variable and encode the push-forward maps.

\paragraph{\bf Contributions.} Building upon the insights from \cite{Bensoussan2020controlHilbert}, we center our approach on Hilbert spaces of push-forward maps. Notably, no prior knowledge in mean-field control is required, as we provide a direct proof of the optimal solution based on the well-known method of completion of square, commonly used in standard control theory for LQ problems. While \cite{Bensoussan2020controlHilbert} considered fully nonlinear objective functions and simple individual linear dynamics, we allow for interaction terms in the dynamics, but we have to restrict the running cost $L$ to be quadratic. Indeed the quadratic structure of $L$ defines a problem-based Sobolev-like norm on $L^2(0,T;L_{\mu^0}^2(\R^n,\R^n))$. This norm is related to the appropriate Hilbertian subspace of controlled push-forward maps that the problem involves. We prove that this subspace is actually a RKHS, with a kernel expressed in closed-form. This property allows us to deal easily with nonlinear objective terms. We also address stochastic dynamics, for which the RKHS functions take their values in random fields, more precisely random push-forward maps. This seems to be the first case where a kernel, a deterministic object, encodes random variables. Finally, through this kernel-based framework, we circumvent the challenges of working in the Wasserstein space of probabilities with second moment; or the difficulty of defining feedback controls in a mean-field context and of finding an adequate Pontryagin's maximum principle. In short, by tackling LQ mean-field control problems through a kernel viewpoint, we want to emphasize that:
\begin{center}
	\emph{Kernels are a powerful formalism which can be applied to many linear-quadratic settings of optimal control, mutatis mutandis.}
\end{center}
The article proceeds gradually in terms of technicality. We present first the deterministic LQ MFC problem in \Cref{sec:lq_mfc} obtaining its solution through a completion of square argument. We then introduce in \Cref{sec:kernel} the corresponding reproducing kernel defining a Hilbert space of time-varying push-forward maps. In \Cref{sec:solution}, we present how to use the kernel to solve LQ MFC problems, also with a final nonlinear term. Finally we show that the kernel formalism extends to stochastic dynamics in the standalone \Cref{sec:stochastic}.

\section{LQ MFC PROBLEMS WRITTEN ON PUSH-FORWARD MAPS}\label{sec:lq_mfc}

\subsection{PRELIMINARIES}
Let $m$ be a probability distribution on $\R^{n}$. We set
\begin{equation}
H:=L_{m}^2(\R^{n},\R^n)=\{X:x\in \R^{n}\mapsto X_{x}\in \R^{n}\,|\,\int_{\R^{n}}|X_{x}|^{2}dm(x)<+\infty\}.\label{eq:2-1}
\end{equation}
The space $H$ is a Hilbert space, which will be the state space. We similarly define a control space 
\begin{equation}
U:=L_{m}^2(\R^{n},\R^{d})=\{V:x\in \R^{n}\mapsto V_{x}\in \R^{d}\,|\,\int_{\R^{n}}|V_{x}|^{2}dm(x)<+\infty\}\label{eq:2-2}
\end{equation}
The norms in $H$ and $U$ will be denoted by $|X|_{H}$ and $|V|_{U}$ respectively. We introduce also the averages 
\begin{equation}
\Xb=\int_{\R^{n}}X_{x}dm(x),\:\Vb=\int_{\R^{n}}V_{x}dm(x).\label{eq:2-3}
\end{equation} 
Consider next the matrices $F(\cdot),\Fb(\cdot)\in L^\infty(t,T;\cL(\R^{n},\R^{n}))$, $G(\cdot)\in L^\infty(t,T;\cL(\R^{d},\R^{n}))$ and vector field $f(\cdot)\in L^2(t,T;\R^{n})$. We consider the linear dynamic system, whose state is in $\R^n$, defined by the equation 
\begin{equation}
\dfrac{d}{ds}X_x(s)=F(s)X_x(s)+\Fb(s)\Xb(s)+G(s)V_{x}(s)+f(s),\, \forall s>t, \quad \quad X_x(t)=X_{x}^t\in\R^n,\, \forall s\in[t,T], x\in\R^n \label{eq:2-4}
\end{equation}
in which $X^t\in H$ is the initial condition and $V(s)=V_{x}(s)\in L^{2}(t,T;U)$ is a control policy. Equation (\ref{eq:2-4}) is easily solved in two steps, by considering separately the evolution of the mean, and that of a given point $X_x^t$. We first notice that $\Xb(s)$ is the solution of 
\begin{equation}
\dfrac{d}{ds}\Xb(s)=(F(s)+\Fb(s))\Xb(s)+G(s)\Vb(s)+f(s),\, \forall s>t, \quad  \quad \Xb(t)=\Xb^t:=\int_{\R^n}X_{x}^t dm(x) \label{eq:2-5}
\end{equation}
and then we can solve (\ref{eq:2-4}), considering $x$ as a fixed parameter and $\Xb(s)$ as a given function. In this way, we do not really solve (\ref{eq:2-4}) as a linear differential equation in $H$, but as a sequence of two ordinary differential equations in $\R^{n}$, one of them being indexed by the parameter $x$. To consider it as a linear differential equation in $H$, we need to introduce some notation. Define the linear operator $\mathcal{F}(s)\in L^{\infty}(t,T;\cL(H;H))$ as follows 
\begin{equation}
(\mathcal{F}(s)X)_{x}=F(s)X_{x}+\Fb(s)\Xb\label{eq:2-6}
\end{equation}
and extend $G(s)\in L^{\infty}(t,T;\cL(\R^{d};\R^{n}))$ to $L^{\infty}(t,T;\cL(U;H))$ by writing 
\begin{equation}
(G(s)V)_{x}=G(s)V_{x}\label{eq:2-7}
\end{equation}
We can then consider the linear differential equation in $H$, as follows
\begin{equation}
\dfrac{d}{ds}X(s)=\mathcal{F}(s)X(s)+G(s)V(s)+f(s), \quad X(t)=X^t.\label{eq:2-8}
\end{equation}

\begin{remark}[Interpretations of $X$]\label{rmk:intepret_X} Maps $X\in H$ are to be understood as moving Dirac masses $\delta_x$ to $\delta_{X_x}$, $X_x$ being a point in $\R^n$, i.e.\ $X$ is a push-forward map. On the other hand, the time index in $X(s)$ indicates that the push-forward takes time $s$ to do so. The probability measure $\mu^0=X^0_\# m$, ``the push-forward of $m$ by $X^0$'', can be seen as the initial condition of the system, and $\mu_s=X(s)_\# m$ as the resulting probability measures at time $s$. For instance if $X^0(x)=a$ is a constant map, then $\mu^0=\delta_a$ and $\mu_s=\delta_{X_a(s)}$. If $X^0(x)=x$ is the identity map, then $\mu^0=m$.
\end{remark}

\begin{remark}[Feedbacks and nonlinearities]
	Let  $v(\cdot,\cdot)$ be a function from $\R^{n}\times(t,T)$ to $\R^{d}$
	and consider the nonlinear equation 
	\begin{equation}
		\dfrac{d}{ds}X_x(s)=F(s)X_x(s)+\Fb(s)\int_{\R^{n}}X_{x'}(s)dm(x')+G(s)v(X_x(s),s)+f(s),\, \forall s>t, \quad  \quad X_x(t)=X_{x}\label{eq:2-71}
	\end{equation}
	We note that, by loosing the linearity in $X_x$ due to $v$, we cannot solve (\ref{eq:2-71}) by the two-step procedure as above. It is a nonlinear integro-differential equation with a parameter $x$. If we can solve it, then we can set 
	\begin{equation}
		V_{x}(s)=v(X_x(s),s)\label{eq:2-72}
	\end{equation}
	so our formalism  $v(x,s)$ includes feedback controls such that we can solve equation (\ref{eq:2-71}). We will see that the optimal control has this property. It is defined by a feedback, for which (\ref{eq:2-71}) can be solved. Note that we could also have considered in \eqref{eq:2-8} that the affine term $f$ depends also nonlinearly on $x$, so long as $f(s)\in H$. Nevertheless fully nonlinear terms with dependence on $X_x(s)$ are out of reach for our linear formalism.
\end{remark}

\subsection{THE LQ MFC PROBLEM }\label{sec:def_control_problem}

We introduce the symmetric matrices $M_{T},\Mb_{T}\in\cL(\R^{n};\R^{n})$, $M(\cdot),\Mb(\cdot)\in L^{\infty}(t,T;\cL(\R^{n};\R^{n}))$, and  symmetric positive definite and invertible $N(\cdot)\in L^{\infty}(t,T;\cL(\R^{d};\R^{d}))$ with $N^{-1}(\cdot)\in L^{\infty}(t,T;\cL(\R^{d};\R^{d}))$. Let also $S(\cdot)\in L^{\infty}(t,T;\cL(\R^{n};\R^{n}))$, $\alpha(\cdot)\in L^{\infty}(t,T;\R^{n}),\beta(\cdot)\in L^{\infty}(t,T;\R^{d})$, $\alpha_T\in \R^{n}$ be affine terms in the payoff
\begin{multline}\label{eq:2-9}
	J_{Xt}(V(\cdot))=\int_{t}^{T}\left[\dfrac{1}{2}\int_{\R^{n}}M(s)X_{x}(s).X_{x}(s)dm(x)+\dfrac{1}{2}(X_{x}(s)-S(s)\Xb(s)).\Mb(s)(X_{x}(s)-S(s)\Xb(s))\right.\\
	\left.+\dfrac{1}{2}\int_{\R^{n}}N(s)V_{x}(s).V_{x}(s)dm(x)	+\alpha(s).\Xb(s)+\beta(s).\Vb(s)\right]ds\\
	+\alpha_T.\Xb(T)+\dfrac{1}{2}\int_{\R^{n}}\left[M_{T}X_{x}(T).X_{x}(T)
	+\dfrac{1}{2}(X_{x}(T)-S(T)\Xb(T)).\Mb(T)(X_{x}(T)-S(T)\Xb(T))\right]dm(x).
\end{multline}
The problem (\ref{eq:2-4}),(\ref{eq:2-9}) is a deterministic linear quadratic mean field control problem. The terms in $M$ and $N$ are classical in LQ optimal control problems, and penalize the deviation from a nominal null trajectory and null control. The extra term in $S$ describes the extra costs incurred if some trajectories originated from $X^0_\# m$ deviate away from the average behavior$\Xb(\cdot)$. We refer to \cite[Section 2]{Bensoussan2015lqmfg} for a more-in-depth discussion of this model and for how it originates from a mean-field limit of a game with finitely-many players. This payoff function was also considered in \cite[Section 6]{Bensoussan2013mfg} and \cite[Section 6.7.1]{Carmona2018}.

We can also express the payoff functional using a control in $L^{2}(t,T;U)$ and state space in $L^{2}(t,T;H)$, with dynamics (\ref{eq:2-8}). Introduce the linear bounded operator $\mathcal{M}(s)\in L^{\infty}(t,T;\cL(H,H))$ defined by 
\begin{equation}
(\mathcal{M}(s)X)_{x}:=(M(s)+\Mb(s))X_{x}+\Mb_S(s)\Xb,\quad \Mb_S(s):=S^{*}(s)\Mb(s)S(s)-S^{*}(s)\Mb(s)-\Mb(s)S(s), \label{eq:2-10}
\end{equation}
and similarly $\mathcal{M}_{T}\in\cL(H,H)$ defined by 
\begin{equation}
(\mathcal{M}_{T}X)_{x}:=(M_{T}+\Mb_{T})X_{x}+\Mb_S(T)\Xb, \quad \Mb_S(T):=S^{*}(T)\Mb_{T}S(T)-S^{*}(T)\Mb_{T}-\Mb_{T}S(T).\label{eq:2-12}
\end{equation}
Then we write $J_{Xt}(V(\cdot))$ as follows 
\begin{multline}
J_{Xt}(V(\cdot))=\int_{t}^{T}\left[\dfrac{1}{2}(\mathcal{M}(s)X(s),X(s))_{H}+\dfrac{1}{2}(N(s)V(s),V(s))_{U}+(\alpha(s),X(s))_{H}+(\beta(s),V(s))_{U}\right]dt\\
+\dfrac{1}{2}(\mathcal{M}_{T}X(T),X(T))_{H}+(\alpha_T,X(T))_{H}\label{eq:2-11}
\end{multline}
We require from now on that the matrices $(M_{T}+\Mb_{T})$,  $M_{T}+\Mb(T)+\Mb_S(T)$, $M(s)+\Mb(s)$ and $(M(s)+\Mb(s)+\Mb_S(s))$, for all $s$, are all positive semidefinite. This is for instance the case if all the matrices $M(\cdot)$ and $\Mb(\cdot)$ are positive semidefinite.

\subsection{SOLUTION OF THE LQ MFC PROBLEM}

We will proceed by completion of square. Since both $X_x(s)$ and $\Xb(s)$ evolve, we introduce two Riccati equations, where $P$ corresponds more or less to $X_x$ and $\Sigma$ to $\Xb$,
\begin{align}
&\dfrac{d}{ds}P(s)+P(s)F(s)+F^{*}(s)P(s)-P(s)G(s)N^{-1}(s)G^{*}(s)P(s)+M(s)+\Mb(s)=0,\quad P(T)=M_{T}+\Mb_{T} \label{eq:2-13} \\
&\dfrac{d}{ds}\Sigma(s)+\Sigma(s)(F(s)+\Fb(s))+(F^{*}(s)+\Fb^{*}(s))\Sigma(s)-\Sigma(s)G(s)N^{-1}(s)G^{*}(s)\Sigma(s)+M(s)+\Mb(s)+\Mb_S(s)=0,\nonumber\\  
&\quad \quad \Sigma(T)=M_{T}+\Mb(T)+\Mb_S(T).\label{eq:2-14}
\end{align}
It is a classical result that the solutions for these equations exist, are unique and are also positive semidefinite. We will also use 
\begin{equation}
\Gamma(s)=\Sigma(s)-P(s)\label{eq:2-15}
\end{equation}
which is the solution of the following Riccati equation
\begin{multline}\label{eq:2-16}
0=\dfrac{d}{ds}\Gamma(s)+\Gamma(s)(F(s)+\Fb(s)-G(s)N^{-1}(s)G^{*}(s)P(s))+(F^{*}(s)+\Fb^{*}(s)-P(s)G(s)N^{-1}(s)G^{*}(s))\Gamma(s)\\
-\Gamma(s)G(s)N^{-1}(s)G^{*}(s)\Gamma(s)
+\Mb_S(s)+P(s)\Fb(s)+\Fb^{*}(s)P(s), \quad \Gamma(T)=\Mb_S(T).
\end{multline}
 We will also need the function $\lambda(s)\in\R^d$ solution of 
\begin{multline}\label{eq:2-17}
0=\dfrac{d}{ds}\lambda(s)+(F^{*}(s)+\Fb^{*}(s)-\Sigma(s)G(s)N^{-1}(s)G^{*}(s))\lambda(s)
+\alpha(s)+\Sigma(s)(f(s)-G(s)N^{-1}(s)\beta(s)), \quad  \lambda(T)=\alpha_T.
\end{multline}
Notice that if $\alpha$, $\beta$, $\alpha_T$ and $f$ are all null, then $\lambda(\cdot)\equiv 0$, so the function $\lambda(\cdot)$ takes care of all the affine terms in $J$. We solve successively 
\begin{gather*}
\dfrac{d}{ds}\hat{\Xb}(s)=(F(s)+\Fb(s)-G(s)N^{-1}(s)G^{*}(s)\Sigma(s))\hat{\Xb}(s)
-G(s)N^{-1}(s)(G^{*}(s)\lambda(s)+\beta(s))+f(s), \quad \hat{\Xb}(t)=\Xb^t\\
\dfrac{d}{ds}\hat{X}_{x}(s)=(F(s)-G(s)N^{-1}(s)G^{*}(s)P(s))\hat{X}_{x}(s)+(F(s)-G(s)N^{-1}(s)G^{*}(s)\Gamma(s))\hat{\Xb}(s)\\-G(s)N^{-1}(s)(G^{*}(s)\lambda(s)+\beta(s))+f(s), \quad \hat{X}_{x}(t)=X_{x}^t,
\end{gather*}
%
to obtain the solution. Indeed we can now state the following theorem
\begin{thm}
\label{theo2-1}The control 
\begin{equation}
\hat{V}_{x}(s)=-N^{-1}(s)G^{*}(s)\left(P(s)\hat{X}_{x}(s)+\Gamma(s)\right)\hat{\Xb}(s))-N^{-1}(s)\left(G^{*}(s)\lambda(s)+\beta(s)\right)\label{eq:2-20}
\end{equation}
is the unique minimum of the functional $J_{Xt}(V(\cdot))$ defined by (\ref{eq:2-9}). Furthermore, the value function satisfies the closed-form expression:
\begin{equation}\label{eq:2-21}
\Phi(X,t):=\inf_{V(\cdot)}J_{Xt}(V(\cdot))=\dfrac{1}{2}\int_{\R^{n}}P(t)X_{x}^t.X_{x}^tdm(x)+\dfrac{1}{2}\Gamma(t)\Xb^t.\Xb^t+\lambda(t).\Xb^t+C_t
\end{equation}
where  $C_t=\int_{t}^{T}f(s).\lambda(s)ds-\dfrac{1}{2}\int_{t}^{T}N^{-1}(s)(G^{*}(s)\lambda(s)+\beta(s)).(G^{*}(s)\lambda(s)+\beta(s))ds$.
\end{thm}

\subsection{PROOF OF THEOREM \ref{theo2-1}}

We define the deviations to the average
\[
\Xt_{x}(s):=X_{x}(s)-\Xb(s),\:\Vt_{x}(s):=V_{x}(s)-\Vb(s),
\]
then we see immediately that 
\begin{equation}
\dfrac{d}{ds}\Xt_{x}(s)=F(s)\Xt_{x}(s)+G(s)\Vt_{x}(s),\quad \Xt_{x}(t)=X_{x}-\Xb=:\Xt_{x}\label{eq:2-22}
\end{equation}
Using that, by definition, $\int\Xt_x(s) dm(x)=0$ and $\int\Vt_x(s) dm(x)=0$, we can write $J_{Xt}(V(\cdot))$ as the sum of two terms, one for the mean and one for the deviation 
\begin{equation}
J_{Xt}(V(\cdot))=\Jt_{\Xt t}(\Vt(\cdot))+\Jb_{\Xb t}(\Vb(\cdot))\label{eq:2-23}
\end{equation}
with
\begin{multline*}
\Jt_{\Xt t}(\Vt(\cdot)):=\int_{\R^{n}}\left\{ \int_{t}^{T}\dfrac{1}{2}\left((M(s)+\Mb(s))\Xt_{x}(s).\Xt_{x}(s)+N(s)\Vt_{x}(s).\Vt_{x}(s)\right)ds+\dfrac{1}{2}(M_{T}+\Mb_{T})\Xt_{x}(T).\Xt_{x}(T)\right\} dm(x)
\end{multline*}
\begin{multline}\label{eq:2-25}
\Jb_{\Xb t}(\Vb(\cdot)):=\int_{t}^{T}\left[\dfrac{1}{2}(M(s)+\Mb(s)+\Mb_S(s))\Xb(s).\Xb(s)+\dfrac{1}{2}N(s)\Vb(s).\Vb(s)+\alpha(s).\Xb(s)+\beta(s).\Vb(s)\right]ds \\
+\dfrac{1}{2}(M_{T}+\Mb_T+\Mb_S(T))\Xb(T).\Xb(T)+\alpha_T.\Xb(T).
\end{multline}
The completion of square works as follows. We define $U_{x}(s)$ as the deviation to our candidate optimal feedback
\begin{equation}
V_{x}(s)=-N^{-1}(s)G^{*}(s)(P(s)X_{x}(s)+\Gamma(s)\Xb(s))-N^{-1}(s)(G^{*}(s)\lambda(s)+\beta(s))+U_{x}(s).\label{eq:2-26}
\end{equation}
It follows that, for $\Ub(s):=\int_{\R^{n}}U_x(s)dm(x)$ and $\Ut_x(s):= U_{x}(s)-\Ub(s)$, using the relation \eqref{eq:2-15},
\begin{align}
\Vb(s)&=-N^{-1}(s)G^{*}(s)\Sigma(s)\Xb(s)-N^{-1}(s)(G^{*}(s)\lambda(s)+\beta(s))+\Ub(s),\label{eq:2-27}\\
\Vt_x(s)&=\Ut_{x}(s)-N^{-1}(s)G^{*}(s)P(s)\Xt_x(s).\nonumber
\end{align}
Consequently
\begin{align*}
	N(s)\Vt_{x}(s).\Vt_{x}(s)&=N(s)\Ut_{x}(s).\Ut_{x}(s)-2\Ut_{x}(s).G^{*}(s)P(s)\Xt_{x}(s)+P(s)G(s) N^{-1}(s)G^{*}(s)P(s)\Xt_{x}(s).\Xt_{x}(s)\\
	&=N(s)\Ut_{x}(s).\Ut_{x}(s)-2G(s)\Vt_x(s).P(s)\Xt_x(s)-P(s)G(s) N^{-1}(s)G^{*}(s)P(s)\Xt_{x}(s).\Xt_{x}(s).
\end{align*}
From \eqref{eq:2-22}, we have that $G(s)\Vt_x(s)=\dfrac{d}{ds}\Xt_{x}(s)-F(s)\Xt_{x}(s)$, hence
\begin{multline*}
	\Jt_{\Xt t}(\Vt(\cdot))=\int_{\R^{n}}\left\{ \int_{t}^{T}\dfrac{1}{2}\left[\left(M(s)+\Mb(s)-P(s)G(s) N^{-1}(s)G^{*}(s)P(s)+F(s)P(s)+P(s)F^*(s)\right)\Xt_{x}(s).\Xt_{x}(s)\right.\right. \\
	\left.\left.-2P(s)\Xt_x(s).\dfrac{d}{ds}\Xt_{x}(s)+N(s)\Ut_{x}(s).\Ut_{x}(s)\right]ds+\dfrac{1}{2}(M_{T}+\Mb_{T})\Xt_{x}(T).\Xt_{x}(T)\right\} dm(x).
\end{multline*}
We recognize $\dfrac{d}{ds}P(s)$ from \eqref{eq:2-13}. We regroup the derivatives into $\dfrac{d}{ds}[P(s)\Xt_{x}(s).\Xt_{x}(s)]$. Integrating, we obtain that 
\begin{equation}
\Jt_{\Xt t}(\Vt(\cdot))=\dfrac{1}{2}\int_{\R^{n}}(X^t_{x}-\Xb^t).P(t)(X^t_{x}-\Xb^t)dm(x)+\dfrac{1}{2}\int_{\R^{n}}\int_{t}^{T}N(s)\Ut_{x}(s).\Ut_{x}(s)ds \,dm(x).\label{eq:2-28}
\end{equation}
Similarly, we recognize in \eqref{eq:2-25} part of $\dfrac{d}{ds}\Sigma(s)$ from \eqref{eq:2-14}. Hence, developing the square terms and hiding the dependence in $s$ for ease of reading,
\begin{multline*}
	\Jb_{\Xb t}(\Vb(\cdot))=\dfrac{-1}{2}\int_{t}^{T}\left[\left(\dfrac{d}{ds}\Sigma+\Sigma(F+\Fb)+(F^{*}+\Fb^{*})\Sigma-\Sigma GN^{-1}G^{*}\Sigma\right)\Xb.\Xb-\Sigma GN^{-1}G^{*}\Sigma\Xb.\Xb\right] ds
	\\
	-C_t+\int_t^Tf.\lambda(s) ds+\int_{t}^{T}\left[\dfrac{1}{2}N\Ub.\Ub-G^{*}\Sigma\Xb.(\Vb+N^{-1}G^{*}\Sigma\Xb)-(G^*\lambda+\beta).\Ub+\alpha.\Xb+\beta.\Vb\right]ds \\
	+\dfrac{1}{2}\Sigma(T)\Xb(T).\Xb(T)+\alpha_T.\Xb(T).
\end{multline*}
Now we notice that $\dfrac{d}{ds}\Xb(s)$ from \eqref{eq:2-5} appears in the first term and that some terms compensate, so
\begin{multline*}
	\Jb_{\Xb t}(\Vb(\cdot))=\dfrac{1}{2}\int_{t}^{T}\left[\dfrac{d}{ds}\left(-\Sigma(s)\Xb(s).\Xb(s)\right)\right] ds+\int_{t}^{T}\left[\dfrac{1}{2}N\Ub.\Ub+\alpha.\Xb+\beta.(\Vb-\Ub)+f.\Sigma\Xb(s)-G^*\lambda.\Ub(s)\right]ds
	\\
	-C_t+\int_t^Tf.\lambda(s) ds +\dfrac{1}{2}\Sigma(T)\Xb(T).\Xb(T)+\alpha_T.\Xb(T).
\end{multline*}
Furthermore, from \eqref{eq:2-17}, we deduce that
\begin{gather*}
	f.\Sigma\Xb+\alpha.\Xb=-\dfrac{d}{ds}(\lambda.\Xb)+(G^*\lambda+\beta).N^{-1}G^*\Sigma\Xb-\lambda.\left(\dfrac{d}{ds}\Xb-G\Vb-f\right)\\
	f.\Sigma\Xb+\alpha.\Xb+\beta.(\Vb-\Ub)-G^*\lambda.\Ub=-\dfrac{d}{ds}(\lambda.\Xb)+\lambda.f+(G^*\lambda+\beta).(\Vb-\Ub+N^{-1}G^*\Sigma\Xb).
\end{gather*}
Finally, assembling the above equations with \eqref{eq:2-27}, we obtain
\begin{equation}
\Jb_{\Xb t}(\Vb(\cdot))=\dfrac{1}{2}\Sigma(t)\Xb^t.\Xb^t+\lambda(t).\Xb^t+\dfrac{1}{2}\int_{t}^{T}N(s)\Ub(s).\Ub(s)ds+C_t\label{eq:2-29}.
\end{equation}
The control $\hat{V}_{x}(\cdot)$ defined by (\ref{eq:2-20}) has a corresponding deviation $\hat{U}_{x}(\cdot)\equiv0$, and is thus clearly the unique minimum for $\Jt_{\Xt t}(\Vt(\cdot))$ and $\Jb_{\Xb t}(\Vb(\cdot))$ based on the obtained expressions \eqref{eq:2-28} and \eqref{eq:2-29}. $\blacksquare$ 

\subsection{ALTERNATIVE EXPRESSION FOR THEOREM \ref{theo2-1}}

Alternatively, we can look at the formulation (\ref{eq:2-8}), (\ref{eq:2-11}) as a control problem in the space $H$. We introduce the Riccati equation with solution $\cP(s)\in\cL(H,H)$, given by 
\begin{equation}
\dfrac{d}{ds}\cP(s)+\cP(s)\mathcal{F}(s)+\mathcal{F^{*}}(s)\cP(s)-\cP(s)G(s)N^{-1}(s)G^{*}(s)\cP(s)+\mathcal{M}(s)=0,\quad \quad \cP(T)=\mathcal{M}_{T}\label{eq:2-31}
\end{equation}
which has a unique solution. We can identify it as 
\begin{equation}
(\cP(s)X)_{x}=P(s)X_{x}+\Gamma(s)\Xb\label{eq:2-32}
\end{equation}
Indeed, equation (\ref{eq:2-31}) is to be understood in its weak form:
\begin{multline}
(\dfrac{d}{ds}\cP(s)X,Y)_{H}+(\mathcal{F}(s)X,\cP(s)Y)_{H}+(\mathcal{F}(s)Y,\cP(s)X)_{H} \\
-(N^{-1}(s)G^{*}(s)\cP(s)X,G^{*}(s)\cP(s)Y)_{U}+(\mathcal{M}(s)X,Y)_{H}=0, \quad \quad 
(\cP(T)X,Y)_{H}=(\mathcal{M}_{T}X,Y)_{H}. \label{eq:2-33}
\end{multline}
To check (\ref{eq:2-32}), we need to verify that $P(s)X_{x}+\Gamma(s)\Xb$ indeed solves \eqref{eq:2-33}, in other words that the following holds:
\begin{multline}\label{eq:2-34}
0=\int_{\R^{n}}\dfrac{d}{ds}P(s)X_{x}.Y_{x}dm(x)+\dfrac{d}{ds}\Gamma(s)\Xb.\Yb+\int_{\R^{n}}(F(s)X_{x}+\Fb(s)\Xb).(P(s)Y_{x}
+\Gamma(s)\Yb)dm(x)\\+\int_{\R^{n}}(F(s)Y_{x}+\Fb(s)\Yb).(P(s)X_{x}+\Gamma(s)\Xb)dm(x)-\int_{\R^{n}}N^{-1}(s)G^{*}(s)(P(s)X_{x}
+\Gamma(s)\Xb).G^{*}(s)(P(s)Y_{x}+\Gamma(s)\Yb)dm(x)\\
+\int_{\R^{n}}(M(s)+\Mb(s))X_{x}.Y_{x}dm(x)+\Mb_S(s)\Xb.\Xb.
\end{multline}
By regrouping the terms in $X$ and $Y$, \eqref{eq:2-34} is verified in view of the Riccati equations for $P$ and $\Gamma$ (\ref{eq:2-13}),(\ref{eq:2-16}). We can then interpret the optimal control (\ref{eq:2-20}) as
\begin{equation}
\hat{V}(s)=-N^{-1}(s)G^{*}(s)\cP(s)\hat{X}(s)-N^{-1}(s)(G^{*}(s)\lambda(s)+\beta(s)).\label{eq:2-35}
\end{equation}
This is just a more compact rewriting of \eqref{eq:2-26}. Performing the same completion of square and computations which led to \eqref{eq:2-28} and \eqref{eq:2-29}, we obtain that the payoff $J_{Xt}(V(\cdot))$, defined by (\ref{eq:2-11}) reads as
\begin{equation}
	J_{Xt}(V(\cdot))=\dfrac{1}{2}(\cP(t)X,X)_{H}+(\lambda(t),X)_{H}+\dfrac{1}{2}\int_{t}^{T}(N(s)U(s),U(s))_{U}ds+C_t,\label{eq:2-37}
\end{equation}
from which we obtain as previously that the control $\hat{V}(s)$ defined by (\ref{eq:2-35}) is the unique optimal control.

\section{KERNEL ASSOCIATED WITH THE LQ MFC PROBLEM}\label{sec:kernel}


The completion of squares has the limitation that it is restricted to quadratic terminal costs. To address nonlinear payoffs, in line with our previous work \cite{aubin2022operator}, we will now introduce our key concept: the reproducing kernel. The latter will allow us to solve a problem of calculus of variations directly on $X(\cdot)$ rather than on $V(\cdot)$.

Owing to the results of \cite{Senkene1973} and Theorem 2.12 in \cite{burbea1984banach} on operator-valued kernels, given a Hilbert space $(\cH,\noridx{\cdot}{\cH})$ of functions of $[t_0,T]$ into a separable Hilbert space $H$, such that the strong topology of $\cH$ is stronger than pointwise convergence, there exists a unique function $K:s,t\in[t_0,T]\mapsto K(s,t)\in\cL(H;H)$, which is called the reproducing $\cL(H;H)$-valued kernel of $\cH$, satisfying
\begin{gather}
	K(\cdot,t)\in\cH, \, \forall t\label{eq:3-100} \\
	\forall\:X(\cdot)\in\cH,\;(X(t),Z)_{H}=(X(\cdot),K(\cdot,t)Z)_\cH,\forall t,Z\in H \label{eq:3-100b} 
\end{gather}
Conversely a kernel $K$ satisfying \eqref{eq:3-100}-\eqref{eq:3-100b} characterizes $\cH$. We then say that $\cH$ is a $H$-valued reproducing kernel Hilbert space (RKHS), and that $K(\cdot,\cdot)$ is the $\cL(H;H)$-valued kernel associated with $\cH$. We refer to \cite{caponetto2008universal} for a more recent presentation of the topic.

We want to show that the Linear Quadratic Mean Field Control problem (\ref{eq:2-8}),(\ref{eq:2-11}) can be associated with a $H$-operator-valued kernel. We consider the subset $\cH$ of $L^{2}(0,T;H)$ defined as follows: 
\begin{equation}
\cH=\{X(\cdot) \, | \, \dfrac{dX}{ds}=\mathcal{F}(s)X(s)+G(s)V(s),\text{with}\;V(\cdot)\in L^{2}(0,T;U)\}\label{eq:3-1}
\end{equation}
equipped with the norm 
\begin{equation}
\|X(\cdot)\|_{\cH}^{2}=(\cJ_0 X(0),X(0))_{H}+(\cJ_T X(T),X(T))_{H}+\int_{0}^{T}(\mathcal{M}(s)X(s),X(s))_{H}dt+\int_{0}^{T}(N(s)V(s),V(s))_{U}ds\label{eq:3-2}
\end{equation} 
where $\cJ_0 $ and $\cJ_T$ are symmetric positive semidefinite linear operators on $H$, and $\cJ_0 $ is assumed to be invertible. For a given $X(\cdot)$ in $\cH$, if $G(s)$ is not injective, we may have several $V(\cdot)$ corresponding to $X(\cdot)$. We will always assign to $X(\cdot)$ its \textbf{representative control}, namely the one which minimizes $(N(s)V(s),V(s))_{U}$ for any s. Any $\Vt(\cdot)$ which corresponds to $X(\cdot)$ is of the form 
\[
\Vt(s)=V(s)+W(s),
\]
with $G(s)W(s)=0$ for all $s$. The representative control is the unique $V(s)$, corresponding to $X(s)$, satisfying the condition 
\begin{equation}
(N(s)V(s), W(s))_{U}=0,\forall W(s)\:\text{such that }G(s)W(s)=0,\forall s.\label{eq:3-200}
\end{equation}


We now define the semigroup on $H$, denoted by $\Phi_{\cF,\cP}(s,t),s>t$,
related to the differential equation 
\begin{equation}
\dfrac{d}{ds}\Phi_{\cF,\cP}(s,t)=(\mathcal{F}(s)-G(s)N^{-1}(s)G^{*}(s)\cP(s))\Phi_{\cF,\cP}(s,t),\quad \Phi_{\cF,\cP}(s,s)=\text{Id}_H\label{eq:3-201}
\end{equation}
where $\cP(s)$ is the solution of the Riccati equation in $H$ generalizing \eqref{eq:2-31}
\begin{equation}
\dfrac{d}{ds}\cP(s)+\cP(s)\mathcal{F}(s)+\mathcal{F^{*}}(s)\cP(s)-\cP(s)G(s)N^{-1}(s)G^{*}(s)\cP(s)+\mathcal{M}(s)=0, \quad \cP(T)=\cJ_T. \label{eq:3-202}
\end{equation}
The main result is the following 
\begin{thm}
\label{theo3-1} Define the family of operators $K(s,t)\in\cL(H,H)$ by the formula 
\begin{equation}
K(s,t):=\Phi_{\cF,\cP}(s,0)(\cJ_0 +\cP(0))^{-1}\Phi_{\cF,\cP}^{*}(t,0)+\int_{0}^{\min(s,t)}\Phi_{\cF,\cP}(s,\tau)G(\tau)N^{-1}(\tau)G^{*}(\tau)\Phi_{\cF,\cP}^{*}(t,\tau)d\tau\label{eq:3-3}
\end{equation}
then the space $\cH$ defined by (\ref{eq:3-1}),(\ref{eq:3-2}) is a RKHS associated with the kernel $K$. 
\end{thm}
\begin{proof} The following proof is similar to that of \cite{aubin2022operator}. Fix $Z\in H$. We first define $\rho_{Zt}(s):=\Phi_{\cF,\cP}^{*}(t,s)Z,s<t$, solution of 
\begin{equation}
-\dfrac{d}{ds}\rho_{Zt}(s)=(\mathcal{F}^{*}(s)-\cP(s)G(s)N^{-1}(s)G^{*}(s))\rho_{Zt}(s)\:s<t,\quad \quad \rho_{Zt}(t)=Z.\label{eq:3-30}
\end{equation}
We next define $\chi_{Zt}(s)$ by 
\begin{equation}
-\dfrac{d}{ds}\chi_{Zt}(s)=\mathcal{F}^{*}(s)\chi_{Zt}(s),\;s<t, \quad \quad \chi_{Zt}(t)=Z.\label{eq:3-31}
\end{equation}
We finally set 
\begin{equation}
r_{Zt}(s):=-(\rho_{Zt}(s)-\chi_{Zt}(s))\1_{s<t}\label{eq:3-32}
\end{equation}
which is defined for any $s\in[0,T]$, and $r_{Zt}(s)=0$ if $s\geq t$. The latter satisfies the following differential equation
\begin{equation}
-\dfrac{d}{ds}r_{Zt}(s)=(\mathcal{F}^{*}(s)-\cP(s)G(s)N^{-1}(s)G^{*}(s))r_{Zt}(s)+\cP(s)G(s)N^{-1}(s)G^{*}(s)\chi_{Zt}(s)\1_{s<t},\quad r_{Zt}(T)=0.\label{eq:3-33}
\end{equation}
We consider next 
\begin{equation}
Y_{Zt}(s):=K(s,t)Z\stackrel{\eqref{eq:3-3}}{=}\Phi_{\cF,\cP}(s,0)(\cJ_0 +\cP(0))^{-1}\rho_{Zt}(0)+\int_{0}^{s}\Phi_{\cF,\cP}(s,\tau)G(\tau)N^{-1}(\tau)G^{*}(\tau)\rho_{Zt}(\tau)\1_{\tau<t}d\tau\label{eq:3-34}
\end{equation}
and from the definition of $\Phi_{\cF,\cP}(s,\tau)$, we can write 
\begin{equation*}
\dfrac{dY_{Zt}(s)}{ds}=(\mathcal{F}(s)-G(s)N^{-1}(s)G^{*}(s)\cP(s))Y_{Zt}(s)+G(s)N^{-1}(s)G^{*}(s)\rho_{zt}(s)\1_{s<t}, \quad Y_{Zt}(0)=(\cJ_0 +\cP(0))^{-1}\rho_{Zt}(0)
\end{equation*}
Using (\ref{eq:3-32}) we can thus write 
\begin{equation}
\dfrac{dY_{Zt}(s)}{ds}=\mathcal{F}(s)Y_{Zt}(s)-G(s)N^{-1}(s)G^{*}(s)(\cP(s)Y_{Zt}(s)+r_{Zt}(s))+G(s)N^{-1}(s)G^{*}(s)\chi_{Zt}(s)\1_{s<t}\label{eq:3-80}
\end{equation}
\[
Y_{Zt}(0)=(\cJ_0 +\cP(0))^{-1}(-r_{Zt}(0)+\chi_{Zt}(0)).
\]
We then introduce the function 
\begin{equation}
Q_{Zt}(s)=\cP(s)Y_{Zt}(s)+r_{Zt}(s)\label{eq:3-11}
\end{equation}
Using the relation (\ref{eq:3-202}) as well as (\ref{eq:3-33}),(\ref{eq:3-80}),
we check easily that the pair $Y_{Zt}(s),Q_{Zt}(s)$ satisfies the
system of forward backward differential equations 
\begin{gather*}
\dfrac{dY_{Zt}(s)}{ds}=\mathcal{F}(s)Y_{Zt}(s)-G(s)N^{-1}(s)G^{*}(s)Q_{Zt}(s)+G(s)N^{-1}(s)G^{*}(s)\chi_{Zt}(s)\1_{s<t},\quad Y_{Zt}(0)=\cJ_0 ^{-1}(\chi_{Zt}(0)-Q_{Zt}(0))\\
-\dfrac{dQ_{Zt}(s)}{ds}=\mathcal{F^{*}}(s)Q_{Zt}(s)+\mathcal{M}(s)Y_{Zt}(s), \quad Q_{Zt}(T)=\cJ_T Y_{Zt}(T). \nonumber
\end{gather*}
Consequently, we see that $Y_{Zt}(s)$ belongs to $\cH$ with corresponding control 
\[
V_{Zt}(s)=-N^{-1}(s)G^{*}(s)(Q_{Zt}(s)-\chi_{Zt}(s)\1_{s<t}).
\]
This control is the representative control of $Y_{Zt}(s)$, since the condition (\ref{eq:3-200}) is clearly satisfied. We can now proceed with the proof that $\cH$ is a RKHS associated to the kernel $K$ defined by (\ref{eq:3-3}). The first property \eqref{eq:3-100} has just been checked. It remains to prove the second property (\ref{eq:3-100b}). If we take $X(\cdot)\in\cH$, see (\ref{eq:3-1}), and denote by $V(\cdot)$ the representative control related to $X(\cdot)$ then we have 
\begin{align*}
	(Y_{Zt}(\cdot),X(\cdot))_\cH&=(X(0),\chi_{Zt}(0)-Q_{Zt}(0))_{H}+\int_{0}^{T}(\mathcal{M}(s)Y_{Zt}(s),X(s))_{X}ds\\
	&\quad-\int_{0}^{T}(N(s)V(s),N^{-1}(s)G^{*}(s)(Q_{Zt}(s)-\chi_{Zt}(s)\1_{s<t}))_{U}ds+(X(T),Q_{Zt}(T))_{H}\\
	&=(X(0),\chi_{Zt}(0)-Q_{Zt}(0))_{H}+\int_{0}^{T}(-\dfrac{dQ_{Zt}(s)}{ds}-\mathcal{F^{*}}(s)Q_{Zt}(s),X(s))_{H}ds\\
	&\quad-\int_{0}^{T}(\dfrac{dX(s)}{ds}-\cF^*(s)X(s),Q_{Zt}(s)-\chi_{Zt}(s)\1_{s<t})_{H}ds +(X(T),Q_{Zt}(T))_{H}\\
	&=(X(0),\chi_{Zt}(0))_{H}+\int_{0}^{t}(G(s)V(s),\chi_{Zt}(s))_{H}ds=(X(t),Z)_{H}
\end{align*}
which completes the proof of (\ref{eq:3-100b}) and of the theorem.
\end{proof}

\section{SOLVING MFC PROBLEMS WITH KERNELS}\label{sec:solution}

\subsection{A RKHS WITH NULL INITIAL CONDITION}

In the definition of $K$, there were two terms. The role of the first was actually to deal with uncontrolled evolutions, so we focus on the second term. From now on, we denote $\cH$ by $\cH_{K}$ and we define
\begin{equation}
K^{1}(s,t):=\int_{0}^{\min(s,t)}\Phi_{\cF,\cP}(s,\tau)G(\tau)N^{-1}(\tau)G^{*}(\tau)\Phi_{\cF,\cP}^{*}(t,\tau)d\tau.\label{eq:3-14}
\end{equation}
Consider the Hilbertian subspace of $\cH_{K}$ of functions with initial value equal to $0$, called $\cH_{K}^{1}$, namely
\begin{equation}
\cH_{K}^{1}=\{X(\cdot) \, | \, \dfrac{dX}{ds}=\mathcal{F}(s)X(s)+G(s)V(s),\;X(0)=0,\text{with}\;V(\cdot)\in L^{2}(0,T;U)\}.\label{eq:3-19}
\end{equation}
The norm of $\cH_{K}^{1}$ is the same as that of $\cH_{K}$ defined in (\ref{eq:3-2}), hence we have
\begin{equation}
\|X(\cdot)\|_{\cH_{K}^{1}}^{2}=(\cJ_T X(T),X(T))_{H}+\int_{0}^{T}(\mathcal{M}(s)X(s),X(s))_{H}dt+\int_{0}^{T}(N(s)V(s),V(s))_{U}ds.\label{eq:3-400}
\end{equation}
It may seem startling to consider a null initial condition in light of Remark~\ref{rmk:intepret_X} which described $X(s)$ as push-forward maps. However in all the following applications $X(\cdot)\in\cH_{K}^{1}$ will be added to an uncontrolled $X_0(\cdot)$ map, which could for instance satisfy $X_0(0)=\Id_{\R^n}$. Hence elements of $\cH_{K}^{1}$ should be understood as fluctuations around the uncontrolled $X_0(\cdot)$. The following result justifies the notation $\cH_{K}^{1}=\cH_{K^{1}}$.
\begin{prop}
\label{prop3-1} The Hilbert space \textup{$\cH_{K}^{1}$} is a RKHS associated with the operator-valued kernel $K^{1}(s,t)$.
\end{prop}
\begin{proof}
The proof is identical to that of Theorem \ref{theo3-1} for $\cH$ and $K$.
\end{proof}

\subsection{KERNEL SOLUTION OF THE LQ MFC PROBLEM}\label{sec:kernel_sol_lq_mfc}

We return to the standard LQ problem (\ref{eq:2-8}), (\ref{eq:2-11}), where the initial state $X^t$ is known. We fix the initial time at $t=0$. We will first put aside all the affine terms. We make a simple change of controls as follows
\begin{equation}
V(s)=\Vt(s)-N^{-1}(s)\beta(s)\label{eq:4-50}
\end{equation}
 and introduce $X_{0}(s)$ solution of 
\begin{equation}
\dfrac{d}{ds}X_{0}(s)=\mathcal{F}(s)X_{0}(s)-G(s)N^{-1}(s)G(s)\beta(s)+f(s),\, X_{0}(0)=X^0\label{eq:4-51}
\end{equation}
 then we have 
\begin{equation}\label{eq:4-52}
X(s)=X_{0}(s)+\xi(s) \quad \text{ with } \quad \dfrac{d}{ds}\xi(s)=\mathcal{F}(s)\xi(s)+G(s)\Vt(s),\quad \xi(0)=0.
\end{equation}
The payoff (\ref{eq:2-11}) with $t=0$ becomes 
\begin{multline}\label{eq:4-8}
J(V(\cdot))=\dfrac{1}{2}\int_{0}^{T}(\mathcal{M}(s)X_{0}(s),X_{0}(s))_{H}ds+\int_{0}^{T}(\alpha(s),X_{0}(s))_{H}ds+\dfrac{1}{2}\int_{0}^{T}(\mathcal{M}(s)\xi(s),\xi(s))_{H}ds\\
+\int_{0}^{T}(\mathcal{M}(s)X_{0}(s)+\alpha(s),\xi(s))_{H}ds+\dfrac{1}{2}\int_{0}^{T}(N(s)\Vt(s),\Vt(s))_{U}ds-\dfrac{1}{2}\int_{0}^{T}(N^{-1}(s)\beta(s),\beta(s))_{U}ds
\\
+\dfrac{1}{2}(\mathcal{M}_{T}X_{0}(T),X_{0}(T))_{H}+(\alpha_T,X_{0}(T))_{H}+\dfrac{1}{2}(\mathcal{M}_{T}\xi(T),\xi(T))_{H}+(\mathcal{M}_{T}X_{0}(T)+\alpha_T,\xi(T))_{H}.
\end{multline}
The problem amounts to minimizing 
\begin{multline}\label{eq:4-9}
\Jt(\Vt(\cdot))=\dfrac{1}{2}\int_{0}^{T}(\mathcal{M}(s)\xi(s),\xi(s))_{H}ds+\dfrac{1}{2}\int_{0}^{T}(N(s)\Vt(s),\Vt(s))_{U}ds\\
\int_{0}^{T}(\mathcal{M}(s)X_{0}(s)+\alpha(s),\xi(s))_{H}ds+\dfrac{1}{2}(\mathcal{M}_{T}\xi(T),\xi(T))_{H}+(\mathcal{M}_{T}X_{0}(T)+\alpha_T,\xi(T))_{H}.
\end{multline}
Since we want to minimize over $(\xi(\cdot),\Vt(\cdot))$, we may as well take $\Vt(\cdot)$ to be the representative control of the trajectory $\xi(\cdot)$ and minimize over $\cH_{K^{1}}$ the functional 
\begin{equation}
\cJ_0(\xi(\cdot))=\dfrac{1}{2}\|\xi(\cdot)\|_{\cH_{K^{1}}}^{2}+\int_{0}^{T}(\mathcal{M}(s)X_{0}(s)+\alpha(s),\xi(s))_{H}ds+(\mathcal{M}_{T}X_{0}(T)+\alpha_T,\xi(T))_{H}\label{eq:4-90}
\end{equation}
Using the representation property of RKHSs, see (\ref{eq:3-300}), the functional $\cJ_0 (\xi(\cdot))$ can also be written as 
\begin{equation}
\cJ_0 (\xi(\cdot))=\dfrac{1}{2}\|\xi(\cdot)\|_{\cH_{K^{1}}}^{2}+\left(\int_{0}^{T}K^{1}(\cdot,s)(\mathcal{M}(s)X_{0}(s)+\alpha(s))ds+K^1(\cdot,T)(\mathcal{M}_{T}X_{0}(T)+\alpha_T),\xi(\cdot)\right)_{\cH_{K^{1}}}.\label{eq:4-10}
\end{equation}
Taking the derivative in $\xi(\cdot)$ of this convex functional,  we can write the Euler condition of optimality in $\cH_{K^{1}}$ 
\begin{equation}
\hat{\xi}(s)+\int_{0}^{T}K^{1}(s,t)(\mathcal{M}(t)X_{0}(t)+\alpha(t))dt+K^{1}(s,T)(\mathcal{M}_{T}X_{0}(T)+\alpha_T)=0\label{eq:4-120}
\end{equation}
Using formula (\ref{eq:3-14}) we obtain 
\begin{equation*}
\hat{\xi}(s)=-\int_{0}^{s}\Phi_{\cF,\cP}(s,\tau)G(\tau)N^{-1}(\tau)G^{*}(\tau)\left[\underbrace{\int_{\tau}^{T}\Phi_{\cF,\cP}^{*}(t,\tau)(\mathcal{M}(t)X_{0}(t)+\alpha(t))dt+\Phi_{\cF,\cP}^{*}(T,\tau)(\mathcal{M}_{T}X_{0}(T)+\alpha_T)}_{=:r(\tau)}\right],
\end{equation*}
then we have 
\begin{gather}
-\dfrac{d}{ds}r(s)=(\mathcal{F}^{*}(s)-\cP(s)G(s)N^{-1}(s)G^{*}(s))r(s)+\mathcal{M}(s)X_{0}(s)+\alpha(s), \quad r(T)=\mathcal{M}_{T}X_{0}(T)+\alpha_T,\label{eq:4-14} \\
\dfrac{d}{ds}\hat{\xi}(s)=(\mathcal{F}(s)-G(s)N^{-1}(s)G^{*}(s)\cP(s))\hat{\xi}(s)-G(s)N^{-1}(s)G^{*}(s)r(s), \quad \hat{\xi}(0)=0.\label{eq:4-16}
\end{gather}
Define 
\begin{equation}
\hat{\chi}(s):=\cP(s)\hat{\xi}(s)+r(s)\label{eq:4-17}
\end{equation}
then we see easily that the pair $\hat{\xi}(s),\hat{\chi}(s)$
is the solution of the forward-backward system of differential equations
in $H$
\begin{gather*}
\dfrac{d}{ds}\hat{\xi}(s)=\mathcal{F}(s)\hat{\xi}(s)-G(s)N^{-1}(s)G^{*}(s)\hat{\chi}(s),\quad \hat{\xi}(0)=0\\
-\dfrac{d}{ds}\hat{\chi}(s)=\mathcal{F}^{*}(s)\hat{\chi}(s)+\mathcal{M}(s)\hat{\xi}(s)+\mathcal{M}(s)X_{0}(s)+\alpha(s),\quad \hat{\chi}(T)=\mathcal{M}_{T}X_{0}(T)+\alpha_T
\end{gather*}
This is exactly the system of necessary and sufficient conditions of optimality for the control problem in $H$, (\ref{eq:4-52}), (\ref{eq:4-9}) and the optimal control is given by 
\begin{equation}
\hat{V}(s)=-N^{-1}(s)G^{*}(s)\hat{\chi}(s).\label{eq:4-19}
\end{equation}
We obtained thus that $\hat{\xi}(s)$ is the optimal state and $\hat{V}(s)$ is its representative control. 

\subsection{LQ MFC PROBLEMS WITH NONLINEAR TERMS}
We now tackle nonlinear terminal costs. Consider again the state equation in $H$
\begin{equation}
\dfrac{d}{ds}X(s)=\mathcal{F}(s)X(s)+G(s)V(s)+f(s), \, X(0)=X^0\label{eq:4-21}
\end{equation}
We use the Wasserstein metric space of probabilities on $\R^{n}$ with finite second moment, denoted $\cP_{2}(\R^{n})$. We associate with $X(s)$ a probability on $\R^{n}$ denoted by $\mu_s$. It is the push forward of $m$ by the map $X(s):x\mapsto X_{x}(s)$, denoted $\mu_s=X(s)_\#m$. We thus have, considering test functions $\varphi(\xi):\R^{n}\mapsto \R$, continuous and bounded
\begin{equation}
\int_{\R^{n}}\varphi(\xi)d\mu_s(\xi)=\int_{\R^{n}}\varphi(X_{x}(s))dm(x)\label{eq:4-22}
\end{equation}
We now consider functionals $m\mapsto\Phi(m)$ on the Wasserstein space. The continuity is naturally defined since the Wasserstein space is a metric space. We use the concept of functional derivative denoted by $\dfrac{d}{dm}\Phi(m)(\xi),\xi\in \R^{n}$ and defined as follows \cite[Definition 5.43]{Carmona2018} 
\begin{equation}
\Phi(m')-\Phi(m)=\int_{0}^{1}\int_{\R^n}\dfrac{d}{dm}\Phi(m+\lambda(m'-m))(\xi)(dm'(\xi)-dm(\xi))d\lambda\label{eq:4-23}
\end{equation}
for any pair $m,m'$ in $\cP_{2}(\R^{n})$. If such a function $m,\xi\mapsto\dfrac{d}{dm}\Phi(m)(\xi)$ exists, then it is called the functional derivative at $m$ of $\Phi(m)$. Note that it is defined up to a constant. Its gradient in $\xi$, denoted by $D\dfrac{d}{dm}\Phi(m)(\xi)$, is uniquely defined if it exists. The important property we are going to use is the following. If $X\in H$ and we consider $X\mapsto\Phi(X_\#m)$ as a functional on $H$, its Gâteaux derivative is $D\dfrac{d}{dm}\Phi(X_\#m)(X_{x})$, provided 
\begin{equation}
\int_{\R{^n}}\left|D\dfrac{d}{dm}\Phi(X_\#m)(X_{x})\right|^{2}dm(x)<+\infty\label{eq:4-24}
\end{equation}
We can now define a mean field control problem with non-quadratic payoff. The state equation is given by (\ref{eq:4-21}) and the payoff is 
\begin{equation*}
J_{Xt}(V(\cdot))=\int_{t}^{T}\left[\dfrac{1}{2}(\mathcal{M}(s)X(s),X(s))_{H}+\dfrac{1}{2}(N(s)V(s),V(s))_{U}+(\alpha(s),X(s))_{H}+(\beta(s),V(s))_{U}\right]dt+\Phi(X(T)_\#m).
\end{equation*} 
Introducing the function $X_{0}(s)$ defined by (\ref{eq:4-51}), and proceeding as for (\ref{eq:4-90}) with $\cJ_T=0$ we obtain the problem of minimizing over $\cH_{K^{1}}$ the functional 
\begin{equation}
\cJ_0 (\xi(\cdot))=\dfrac{1}{2}\|\xi(\cdot)\|_{\cH_{K^{1}}}^{2}+\int_{0}^{T}(\mathcal{M}(s)X_{0}(s)+\alpha(s),\xi(s))_{H}ds+\Phi((X_{0}(T)+\xi(T))_\#m)\label{eq:4-26}
\end{equation}
So we have also 
\begin{equation}
\cJ_0 (\xi(\cdot))=\dfrac{1}{2}\|\xi(\cdot)\|_{\cH_{K^{1}}}^{2}+\left(\int_{0}^{T}K^{1}(\cdot,s)(\mathcal{M}(s)X_{0}(s)+\alpha(s))ds,\xi(\cdot)\right)_{\cH_{K^{1}}}+\Phi((X_{0}(T)+\xi(T))_\#m)\label{eq:4-27}
\end{equation}
We can write the Euler equation of optimality, i.e.\ for all $\xi(\cdot)\in\cH_{K^{1}}$,
\begin{equation*}
(\hat{\xi}(\cdot),\xi(\cdot))_{\cH_{K^{1}}}+(\int_{0}^{T}K^{1}(\cdot,s)(\mathcal{M}(s)X_{0}(s)+\alpha(s))ds,\xi(\cdot))_{\cH_{K^{1}}}+(K^{1}(\cdot,T)\,D\dfrac{d}{dm}\Phi(\hat{m}(T))(X_{0}(T)+\hat{\xi}(T)),\xi(\cdot))_{\cH_{K^{1}}}=0 
\end{equation*}
in which we have used the notation 
\begin{equation}
\hat{m}(T):=(X_{0}(T)+\hat{\xi}(T))_\#m\label{eq:4-29}
\end{equation}
Therefore $\hat{\xi}(\cdot)$ is solution of the equation 
\begin{equation}
\hat{\xi}(s)+\int_{0}^{T}K^{1}(s,t)(\mathcal{M}(t)X_{0}(t)+\alpha(t))dt+K^{1}(s,T)\,D\dfrac{d}{dm}\Phi(\hat{m}(T))(X_{0}(T)+\hat{\xi}(T))=0 \label{eq:4-30}
\end{equation}
So we define $\hat{\xi}(T)$ by the fixed point equation (in $H)$ 
\begin{equation}
\hat{\xi}(T)+\int_{0}^{T}K^{1}(T,t)(\mathcal{M}(t)X_{0}(t)+\alpha(t))dt+D\dfrac{d}{dm}\Phi(\hat{m}(T))(X_{0}(T)+\hat{\xi}(T))=0\label{eq:4-31}
\end{equation}
\begin{remark}[Example of nonlinearities] \label{rmk:nonlinear}
In variational inference, the ubiquitous functional is the relative entropy (a.k.a.\ the Kullback--Leibler divergence), defined as
	\begin{equation}\label{eq:kl_def}
		\KL(\mu|\pi)=\int\ln\left(\dfrac{d\mu}{d\pi}(x)\right)d\mu(x)
	\end{equation}
for $\mu$ absolutely continuous w.r.t.\ $\pi$ and $+\infty$ otherwise. Through Jensen's inequality, it is easy to show that $\KL(\mu|\pi)\ge 0$ and vanishes if and only if $\mu=\pi$. Fix a probability measure $\pi$ over $\R^n$ with a density w.r.t.\ Lebesgue measure. Then it is easy to show that
\begin{itemize}
	\item for $\Phi(m)=\KL(m | \pi)$, then $\dfrac{d}{dm}\Phi(m)(x)=1+\ln \dfrac{dm}{d\pi}(x)$ for all $m$, where $\dfrac{dm}{d\pi}$ is the Radon-Nikodym derivative, and $+\infty$ if $m$ is not absolutely continuous w.r.t.\ $\pi$. The same holds for the cross-entropy $\Phi(m)=-\int\ln\left(\pi(x)\right)dm(x)$.
	\item for $\Phi(m)=\KL(\pi | m)$, then $\dfrac{d}{dm}\Phi(m)(x)=-\dfrac{d\pi}{dm}(x)$ for all $m$ such that $ \pi$ has a density w.r.t.\ $m$.
\end{itemize}
In these two cases, \eqref{eq:4-30} is a nonlinear equation and $Phi$ is coercive, but not defined everywhere, which requires some care.

Note that we could also have considered in \eqref{eq:4-51} and \eqref{eq:4-8} that the affine terms $f$, $\beta$, $\mu$, $\alpha_T$ depend also nonlinearly on $x$, so for instance on $X_0(s)(x)$. However the Hilbertian structure is quite rigid, and it requires all the other terms to be linear or quadratic in order to be absorbed into the definition of $\cH_1$.\\

We can obviously recover the quadratic case that we dealt with in \Cref{sec:kernel_sol_lq_mfc} if we take 
\begin{equation}
\Phi(m)=\dfrac{1}{2}\int_{\R^{n}}(M_{T}+\Mb_{T})x.xdm(x)+\dfrac{1}{2}\Mb_S(T)\int_{\R^{n}}xdm(x).\int_{\R^{n}}xdm(x)+\alpha_T.\int_{\R^{n}}xdm(x)\label{eq:4-32}
\end{equation}
We then have 
\begin{equation*}
\dfrac{d}{dm}\Phi(m)(x)=\dfrac{1}{2}(M_{T}+\Mb_{T})x.x+\Mb_S(T)\left(\int_{\R^{n}}\xi dm(\xi)\right).x+\alpha_T.x
\end{equation*}
Hence 
\begin{gather*}
D\dfrac{d}{dm}\Phi(m)(x)=(M_{T}+\Mb_{T})x+\Mb_S(T)\int_{\R^{n}}\xi dm(\xi)+\alpha_T, \\
D\dfrac{d}{dm}\Phi(\hat{m}(T))(X_{0}(T)+\hat{\xi}(T))=\mathcal{M}_{T}X_{0}(T)+\alpha_T.
\end{gather*}
Then (\ref{eq:4-31}) coincides with (\ref{eq:4-120}) (for $s=T)$. In other words, terminal quadratic costs can be included or not in the definition of the kernels $K$ and $K^1$.
\end{remark}

\section{STOCHASTIC EXTENSION}\label{sec:stochastic}

In this section, we show that the kernel approach remains valid and fruitful even for stochastic dynamics. This section is naturally an extension of what precedes when the Brownian noise is taken to be null. Even though many of the proofs are similar to the previous ones and we repeat only the key arguments, this section can still be read mostly independently from the rest of the article. Importantly, the adaptivity issue of stochastic dynamics requires a nontrivial modification of the kernel. Furthermore, the question of adapting the kernel framework to settings with a controlled diffusion term remains open. Such an extension would allow to deal with Mean-Variance Portfolio Selection problem in this full generality as presented in \cite[Section 5]{Andersson2010}.

\subsection{DEFINITIONS }

We consider an atomless probability space $(\Omega,\mathcal{A},P)$ sufficiently large to contain a standard Wiener process in $\R^n$, denoted $w(t)$, and rich enough to support random variables that are independent of the entire Wiener process. We now introduce two Hilbert spaces of square integrable random variables, both depending on the fixed probability measure $m$
\begin{equation}
H:=L^{2}(\Omega,\mathcal{A},P;L_{m}^{2}(\R^{n};\R^{n})),\;U:=L^{2}(\Omega,\mathcal{A},P;L_{m}^{2}(\R^{n};\R^{d}))\label{eq:5-1}
\end{equation}
Elements of $H$ (resp.\ $U$) are denoted by $X=X_{x}$ (resp.\ $V=V_{x}$), where $X_{x}=X_{x}(\omega)$ is a random variable and 
\begin{equation}
|X|_{H}^{2}=\bE\int_{\R^{n}}|X_{x}|^{2}dm(x)<+\infty, \quad \quad |V|_{U}^{2}=\bE\int_{\R^{n}}|V_{x}|^{2}dm(x)<+\infty\label{eq:5-2}
\end{equation}
Denote the mean of $X$ by $\Xb$, with 
\begin{equation}
\Xb=\bE\int_{\R^{n}}X_{x}dm(x)\label{eq:5-4}
\end{equation}
and similarly for $\Vb$. For all $s\in(0,T)$, we introduce next operators $\mathcal{F}(s)\in\cL(H;H)$ and $G(s)\in\cL(\R^{d};\R^{n})$, extended as an operator in $\cL(U;H)$, as follows 
\begin{equation}
(\mathcal{F}(s)X)_{x}=F(s)X_{x}+\Fb(s)\Xb, \quad \quad (G(s)V)_{x}=G(s)V_{x},\label{eq:5-5}
\end{equation}
where $F(s),\Fb(s)$ are deterministic bounded matrices in $\cL(\R^{n};\R^{n})$. We will use the notation
\begin{equation}
\Xt=X-\Xb, \quad \Vt=V-\Vb. \label{eq:5-7}
\end{equation}
We finally introduce the operators $\mathcal{M}(s),\mathcal{M}_{T}\in\cL(H;H)$
defined by the formulas
\begin{equation}
(\mathcal{M}(s)X)_{x}=(M(s)+\Mb(s))\Xt_{x}+(M(s)+(I-S^{*}(s))\Mb(s)(I-S(s)))\Xb\label{eq:5-8}
\end{equation}
\[
(\mathcal{M}_{T}X)_{x}=(M_{T}+\Mb_{T})\Xt_{x}+(M_{T}+(I-S^{*}(T))\Mb_{T}(I-S(T)))\Xb
\]
where $M(s),\Mb(s),M_{T},\Mb_{T}$ are symmetric bounded (in $s$) matrices in $\cL(H;H)$ and $S(s),S_{T}$ are bounded (in $s$) matrices in $\cL(H;H)$. We require furthermore that the matrices $(M_{T}+\Mb_{T})$,  $(M_{T}+(I-S^{*}(T))\Mb_{T}(I-S(T)))$, $M(s)+\Mb(s)$ and $(M(s)+(I-S^{*}(s))\Mb(s)(I-S(s)))$, for all $s$, are all positive semidefinite. 

\subsection{CONTROL PROBLEM }

We consider the family of $\sigma$-algebras 
\begin{equation}
\cW_{t}^{s}=\sigma(w(\tau)-w(t),t\leq\tau\leq s),t<s\label{eq:5-9}
\end{equation}
inducing the filtration $\cW_{t}=\{\cW_{t}^{s}\}_{s\ge t}$. From now on $X^0$ is an element of $H$, such that, for any $x$, $X^0_{x}$ is independent of the filtration $\cW_{t}$, and we introduce the $\sigma-$algebras 
\begin{equation}
\cW_{Xt}^{s}=\sigma(X^0)\cup\cW_{t}^{s}\label{eq:5-10}
\end{equation} 
generating the filtration $\cW_{Xt}$. We consider the Hilbert spaces $L^{2}(t,T;H)$, $L^{2}(t,T;U)$ and the Hilbertian subspaces $L_{\cW_{Xt}}^{2}(t,T;H),L_{\cW_{Xt}}^{2}(t,T;U)$ of processes adapted to the filtration $\cW_{Xt}$. Let $f(s)\in L^{2}(t,T;\R^{n})$ and $\eta\in\cL(\R^{n};\R^{n})$. If $V(s)\in L_{\cW_{Xt}}^{2}(t,T;U)$, called the control process, we define the state of a dynamic system as $X(s)$ solution of 
\begin{equation}
X(s)=X^0+\int_{t}^{s}(\mathcal{F}(\tau)X(\tau)+G(\tau)V(\tau)+f(\tau))d\tau+\eta(w(s)-w(t))\label{eq:5-11}
\end{equation}
for a given initial condition $X^0\in H$. We have $X(s)\in L_{\cW_{Xt}}^{2}(t,T;H)$. We then define the payoff 
\begin{multline}
J_{Xt}(V(\cdot)):=\int_{t}^{T}\left[\dfrac{1}{2}(\:(\mathcal{M}(s)X(s),X(s))_H +(N(s)V(s),V(s)))_U +(\alpha(s),X(s))_H +(\beta(s),V(s))_U \right]ds\\
+\dfrac{1}{2}(\mathcal{M}_{T}X(T),X(T))_H +(\alpha_T,X(T))_H \label{eq:5-12}
\end{multline}
 where $N(s)$ is a symmetric positive invertible matrix and $\alpha(s)\in L^{2}(t,T;\R^{n}),\beta(s)\in L^{2}(t,T;\R^{d})$$,\alpha_T\in \R^{n}$.
The objective is to minimize $J_{Xt}(V(\cdot))$ on $L_{\cW_{Xt}}^{2}(t,T;U)$. 

\subsection{SOLUTION OF THE CONTROL PROBLEM }

As done previously in \eqref{eq:2-31}, we introduce a Riccati equation in $H$,
\begin{equation}
\dfrac{d}{ds}\cP(s)+\cP(s)\mathcal{F}(s)+\mathcal{F^{*}}(s)\cP(s)-\cP(s)G(s)N^{-1}(s)G^{*}(s)\cP(s)+\mathcal{M}(s)=0, \, \cP(T)=\mathcal{M}_{T}\label{eq:5-13}
\end{equation}
 whose solution $\cP(s)$ is a symmetric positive linear operator
in $\cL(H,H)$. We also define $\lambda(s)\in L^{2}(t,T;\R^{n})$, the unique
solution of the linear backward equation 
\begin{equation}
-\dfrac{d\lambda}{ds}(s)=\mathcal{F^{*}}(s)\lambda(s)-\cP(s)G(s)N^{-1}(s)(G^{*}(s)\lambda(s)+\beta(s))+\alpha(s)+\cP(s)f(s), \quad \lambda(T)=\alpha_T.\label{eq:5-14}
\end{equation}
A priori $\lambda(s)\in L^{2}(t,T;H)$. However, it is readily seen that $\lambda(s)=\overline{\lambda}(s)$ since $\alpha(s),\alpha_T\in \R^n$ and $\beta(s)\in\R^d$, so integrating \eqref{eq:5-14} over $m$ entails that $\overline{\lambda}$ is also a solution of \eqref{eq:5-14}. Therefore $\lambda(s)$ is simply in $L^{2}(t,T;\R^{n})$. We next define $\hat X(s)\in L_{\cW_{Xt}}^{2}(t,T;H)$, by the equation 
\begin{multline}
\hat X(s)=X^0+\int_{t}^{s}(\mathcal{F}(\tau)-G(\tau)N^{-1}(\tau)G^{*}(\tau)\cP(\tau))\hat X\tau))d\tau\\ +\int_{t}^{s}(f(\tau)-G(\tau)N^{-1}(\tau)(G^{*}(\tau)\lambda(\tau)+\beta(\tau))d\tau+\eta(w(s)-w(t))\label{eq:5-15}
\end{multline}
and set 
\begin{equation}
\hat V(s):=-N^{-1}(s)G^{*}(s)\cP(s)\hat X(s)-N^{-1}(s)(G^{*}(s)\lambda(s)+\beta(s))\label{eq:5-16}
\end{equation}
which belongs to $L_{\cW_{Xt}}^{2}(t,T;U)$. We can state 
\begin{thm}
\label{theo5-1} The process $\hat V(s)$ is the optimal control 
\begin{equation}
J_{Xt}(\hat V(\cdot))=\inf_{V(\cdot)}J_{Xt}(V(\cdot))\label{eq:5-17}
\end{equation}
\end{thm}

\begin{proof} 
This proof using the method of completion of square is similar to the one of Theorem~\ref{theo2-1} up to the stochastic terms. For any control $V(\cdot)$ we define the mean $\Vb(\cdot)$ and deviation $\Vt(\cdot)$=$V(\cdot)-\Vb(\cdot)$.
We can then write 
\begin{equation}
J_{Xt}(V(\cdot))=J_{\Xt t}(\Vt(\cdot))+J_{\Xb t}(\Vb(\cdot))\label{eq:5-18}
\end{equation}
 with 
\begin{equation}
\Xt(s)=\Xt_0+\int_{t}^{s}(\mathcal{F}(\tau)\Xt(\tau)+G(\tau)\Vt(\tau))d\tau+\eta(w(s)-w(t))\label{eq:5-19}
\end{equation}
\begin{equation}
\dfrac{d}{ds}\Xb(s)=\mathcal{F}(s)\Xb(s)+G(s)\Vb(s)+f(s), \quad \Xb(t)=\Xb^0\label{eq:5-20}
\end{equation}
and 
\begin{equation}
J_{\Xt t}(\Vt(\cdot))=\dfrac{1}{2}\int_{t}^{T}\left[(\mathcal{M}(s)\Xt(s),\Xt(s))_H +(N(s)\Vt(s),\Vt(s))_U \right]ds+\dfrac{1}{2}(\mathcal{M}_{T}\Xt(T),\Xt(T))_H \label{eq:5-21}
\end{equation}
\begin{multline}
J_{\Xb t}(\Vb(\cdot))=\int_{t}^{T}\left[\dfrac{1}{2}(\:(\mathcal{M}(s)\Xb(s),\Xb(s))_H +(N(s)\Vb(s),\Vb(s)))_U +(\alpha(s),\Xb(s))_H +(\beta(s)\Vb(s),\Vb(s))_U \right]ds\\
 +\dfrac{1}{2}(\mathcal{M}_{T}\Xb(T),\Xb(T))_H +(\alpha_T,\Xb(T))_H \label{eq:5-22}
\end{multline}
We perform the change of control in (\ref{eq:5-19}) 
\begin{equation}
\Vt(s)=-N^{-1}(s)G^{*}(s)\cP(s)\Xt(s)+\Vt_1(s)\label{eq:5-23}
\end{equation}
Performing the change in (\ref{eq:5-21}), the same computations as in the proof of Theorem~\ref{theo2-1}, except for the Itô integration, lead to
\begin{equation}
J_{\Xt t}(\Vt(\cdot))=\dfrac{1}{2}(\cP(t)\Xt,\Xt)_H +\dfrac{1}{2}\int_{t}^{T}\sum_{j}\cP(s)\eta^{j}.\eta^{j}ds+\dfrac{1}{2}\int_{t}^{T}(N(s)\Vt_1(s),\Vt_1(s))_U ds\label{eq:5-24}
\end{equation}
where $\eta^{j}$ are the column vectors of the matrix $\eta$. Similarly, for equation (\ref{eq:5-22}) we perform the change of control 
\begin{equation}
\Vb(s)=-N^{-1}(s)G^{*}(s)\cP(s)\Xb(s)-N^{-1}(s)(G^{*}(s)\lambda(s)+\beta(s))+\Vb_{1}(s)\label{eq:5-25}
\end{equation}
We obtain, performing the change in (\ref{eq:5-22}) and recalling that $C_t$ was defined in \eqref{eq:2-21},  
\begin{equation}
J_{\Xb t}(\Vb(\cdot))=\dfrac{1}{2}\cP(t)\Xb^0.\Xb^0+\lambda(t).\Xb^0+\dfrac{1}{2}\int_{t}^{T}N(s)\Vb_{1}(s).\Vb_{1}(s)ds+C_t\label{eq:5-26}
\end{equation}
Combining (\ref{eq:5-24}) and (\ref{eq:5-26}) we obtain
\begin{equation}
J_{Xt}(V(\cdot))=\dfrac{1}{2}\cP(t)X^0.X^0+\lambda(t).X^0+C_t+\dfrac{1}{2}\int_{t}^{T}\sum_{j}\cP(s)\eta^{j}.\eta^{j}ds
+\dfrac{1}{2}\int_{t}^{T}(N(s)V_1(s),V_1(s))_Uds\label{eq:5-27}
\end{equation}
and clearly $\hat V(s)$ is optimum since it has corresponding $V_1(\cdot)\equiv 0$. This concludes the proof. 

\end{proof}


\subsection{A RKHS OF RANDOM FIELD-VALUED FUNCTIONS}

From now on, we take the origin $t=0$. The filtration $\cW_{X0}$ is simply denoted by $\cW_{X}$. We will use the spaces $L_{\cW_{X}}^{2}(0,T;H),L_{\cW_{X}}^{2}(0,T;U)$. We introduce the subspace of $L_{\cW_{X}}^{2}(0,T;H)$ defined by
\begin{equation}
	\cH_1=\left\{ \xi(\cdot)\in L_{\cW_{X}}^{2}(0,T;H)\, | \,  \dfrac{d}{ds}\xi(s)=  \mathcal{F}(s)\xi(s)+G(s)v(s), \, \xi(0)=0,\, v(\cdot)\in L_{\cW_{X}}^{2}(0,T;U)\right\} \label{eq:6-1}
\end{equation}
To deal with the situation where there are multiple $v(\cdot)$ for the same $\xi(\cdot)$, we define the \textbf{representative control}. It minimizes the norm $\int_{0}^{T}(N(s)V(s),V(s))_{U}ds$ for all controls which satisfy the differential equation (\ref{eq:6-1}) for the same $\xi(\cdot)$. In other words it is the unique one which satisfies 
\begin{equation}
(N(s)v(s),\Vt(s))_U =0,\text{if}\:G(s)\Vt(s)=0,\text{a.e. }s\in(0,T)\label{eq:6-2}
\end{equation}
We then equip $\cH_1$ with the norm 
\begin{equation}
\|\xi(\cdot)\|_{\cH_1}^{2}=\int_{0}^{T}[(\mathcal{M}(s)\xi(s),\xi(s))_H +(N(s)v(s),v(s))_U ]ds+(\mathcal{M}_{T}\xi(T),\xi(T))_H \label{eq:6-3}
\end{equation}
We introduce the family of Hilbertian subspaces of $H$, called $H_{\cW_{X}^{t},}$ indexed by $t$, each $H_{\cW_{X}^{t}}$ being the subset of $H$ of random variables $z^{t}$ which are $\cW_{X}^{t}$ measurable. A kernel $K_1(s,t)$ mapping $H_{\cW_{X}^{t}}$ to $H_{\cW_{X}^{s}}$ with the properties
\begin{gather}
K_1(\cdot,t)z_{t}\in\cH_1,\forall t,z_{t}\in H_{\cW_{X}^{t}} \label{eq:6-4} \\
(K_1(\cdot,t)z_{t},\xi(\cdot))_{\cH_1}=(\xi(t),z_{t}),\forall\xi(\cdot)\in\cH_1, \,z_{t}\in H_{\cW_{X}^{t}} \label{eq:6-5}
\end{gather}
is called a reproducing kernel for $\cH_1$. Denoting by $\Phi_{\cF,\cP}(s,\tau),s>\tau$ the semigroup associated with the operator $\mathcal{F}(s)-G(s)N^{-1}(s)G^{*}(s)\cP(s)$, we define the family $K_1(s,t)$ as follows. For $z_{t}\in H_{\cW_{X}^{t}}$, we set
\begin{equation}
K_1(s,t)z_{t}:=\int_{0}^{\min(s,t)}\Phi_{\cF,\cP}(s,\tau)G(\tau)N^{-1}(\tau)G^{*}(\tau)\Phi_{\cF,\cP}^{*}(t,\tau)\bE[z_{t}|\cW_{X}^{\tau}]d\tau\label{eq:6-6}
\end{equation}
Recall that $H=L^{2}(\Omega,\mathcal{A},P;L_{m}^{2}(\R^{n};\R^{n}))$, so $\xi(s)$ or $K(s,t)z_t$ are random fields, even though the values of the kernel are deterministic.
\begin{thm}
\label{theo7-1}The kernel $K_1(s,t)$ is a reproducing kernel for the
Hilbert space $\cH_1$.
\end{thm}

\begin{proof}We set
\[
Y_{z_{t}t}(s):=K_1(s,t)z_{t},\quad \rho_{z_{t}t}(s):=\Phi_{\cF,\cP}^{*}(t,s)\bE[z_{t}|\cW_{X}^{s}]
\]
 and want to show that $Y_{z_{t}t}(\cdot)$$\in\cH_1$. We can write
\[
Y_{z_{t}t}(s)=\int_{0}^{s}\Phi_{\cF,\cP}(s,\tau)G(\tau)N^{-1}(\tau)G^{*}(\tau)\rho_{z_{t}t}(\tau)\1_{\tau<t}d\tau
\]
 which means 
\begin{equation}
\dfrac{d}{ds}Y_{z_{t}t}(s)=(\mathcal{F}(s)-G(s)N^{-1}(s)G^{*}(s)\cP(s))Y_{z_{t}t}(s)+G(s)N^{-1}(s)G^{*}(s)\rho_{z_{t}t}(s)\1_{s<t}, \quad Y_{z_{t}t}(0)=0. \label{eq:6-7}
\end{equation}
This implies that $Y_{z_{t}t}(\cdot)\in\cH_1$, with representative control 
\begin{equation}
V_{z_{t}t}(s)=-N^{-1}(s)G^{*}(s)(\cP(s))Y_{z_{t}t}(s)-\rho_{z_{t}t}(s)\1_{s<t}\label{eq:6-8}
\end{equation}
In addition, according to the representation of martingales with respect to Wiener processes, and the independence of $X^0$ and the Wiener process, we can write 
\[
\bE[z_{t}|\cW_{X}^{s}]=\bE[z_{t}|\sigma(X^0)]+\int_{0}^{s}L_{z_{t}t}(\tau)dw(\tau)
\]
where $L_{z_{t}t}(\cdot)\in L_{\cW_{X}}^{2}(0,T;H)$. Therefore, for $s<t$, $\rho_{z_{t}t}(s)$ satisfies the backward stochastic differential equation 
\begin{equation}
-d\rho_{z_{t}t}(s)=(\mathcal{F}^{*}(s)-\cP(s)G(s)N^{-1}(s)G^{*}(s))\rho_{z_{t}t}(s)ds-\Phi_{\cF,\cP}^{*}(t,s)L_{z_{t}t}(s)dw(s),\:s<t, \quad \rho_{z_{t}t}(t)=z_{t}.\label{eq:6-9}
\end{equation}
We next check the reproducing property. Let $\xi(\cdot)\in\cH_1$. We have 
\[
\dfrac{d}{ds}\xi(s)=\mathcal{F}(s)\xi(s)+G(s)v(s), \quad  \xi(0)=0.
\]
From the definition (\ref{eq:6-3}) and \eqref{eq:6-8} we have 
\begin{multline}
(Y_{z_{t}t}(\cdot),\xi(\cdot))_{\cH_1}=\int_{0}^{T}(\xi(s),\mathcal{M}(s)Y_{z_{t}t}(s))_H ds+(\xi(T),\mathcal{M}_{T}Y_{z_{t}t}(T))_H \\
-\int_{0}^{T}(G(s)v(s),\cP(s)Y_{z_{t}t}(s))_H ds+\int_{0}^{t}(G(s)v(s),\rho_{z_{t}t}(s))_H ds\label{eq:6-11}
\end{multline}
We set $Q_{z_{t}t}(s):=\cP(s)Y_{z_{t}t}(s)$, then we obtain 
\begin{equation}
-\dfrac{d}{ds}Q_{z_{t}t}(s)-\mathcal{F}^{*}(s)Q_{z_{t}t}(s)=\mathcal{M}(s)Y_{z_{t}t}(s)-\cP(s)G(s)N^{-1}(s)G^{*}(s)\rho_{z_{t}t}(s)\1_{s<t}, \quad Q_{z_{t}t}(T)=\mathcal{M}_{T}Y_{z_{t}t}(T)\label{eq:6-12}
\end{equation}
Hence, multiplying \eqref{eq:6-12} by $\xi(s)$ and integrating by parts, we obtain that
\begin{multline}
\int_{0}^{T}(\xi(s),\mathcal{M}(s)Y_{z_{t}t}(s))_Hds+(\xi(T),\mathcal{M}_{T}Y_{z_{t}t}(T))_H-\int_{0}^{T}(G(s)v(s),\cP(s)Y_{z_{t}t}(s))_Hds\\
=\int_{0}^{t}(\xi(s),\cP(s)G(s)N^{-1}(s)G^{*}(s)\rho_{z_{t}t}(s))_H ds \label{eq:6-13}
\end{multline}
Similarly, noticing that there is no Itô integral when integrating by parts because the stochastic evolution of $\xi(s)$ does not include a Brownian motion,
\begin{equation}
(\xi(t),z_{t})_H=\bE\left[\int_{0}^{t} \dfrac{d}{ds}(\xi(s).\rho_{z_{t}t}(s))_Hds\right] =\int_{0}^{t}(G(s)v(s),\rho_{z_{t}t}(s))_Hds+\int_{0}^{t}(\xi(s),\cP(s)G(s)N^{-1}(s)G^{*}(s)\rho_{z_{t}t}(s))_H ds.\label{eq:6-14}
\end{equation}
Combining (\ref{eq:6-13}) and (\ref{eq:6-14}) yields 
\begin{multline*}
	\int_{0}^{T}(\xi(s),\mathcal{M}(s)Y_{z_{t}t}(s))_H ds+(\xi(T),\mathcal{M}_{T}Y_{z_{t}t}(T))_H-\int_{0}^{T}(G(s)v(s),\cP(s)Y_{z_{t}t}(s))_Hds\\
	=(\xi(t),z_{t})_H-\int_{0}^{t}(G(s)v(s),\rho_{z_{t}t}(s))_H ds.
\end{multline*}
and from (\ref{eq:6-11}) we get the reproducing property \eqref{eq:6-5}. This
concludes the proof.
\end{proof}


\subsection{DEFINING THE CONTROL PROBLEM }

We turn back to the control problem (\ref{eq:5-11}),(\ref{eq:5-12}), with initial condition at $t=0$ which reads
\begin{gather}
X(s)=X^0+\int_{0}^{s}(\mathcal{F}(\tau)X(\tau)+G(\tau)V(\tau)+f(\tau))d\tau+\eta w(s)\label{eq:7-1}\\
J(V(\cdot))=\int_{0}^{T}\left[\dfrac{1}{2}(\:(\mathcal{M}(s)X(s),X(s))_H +(N(s)V(s),V(s)))_U +(\alpha(s),X(s))_H +(\beta(s),V(s))_U \right]ds \nonumber\\
+\dfrac{1}{2}(\mathcal{M}_{T}X(T),X(T))+(\alpha_T,X(T))_H \label{eq:7-2}
\end{gather}
The control $V(\cdot)$ belongs to $L_{\cW_{X}}^{2}(0,T;U)$. We define successively $v(s)\in L_{\cW_{X}}^{2}(0,T;U)$, $X_{0}(s)\in L_{\cW_{X}}^{2}(0,T;H)$, by 
\begin{gather}
v(s):=V(s)+N^{-1}(s)\beta(s)\label{eq:7-3} \\
X_{0}(s):=X^0+\int_{0}^{s}(\mathcal{F}(\tau)X_{0}(\tau)+f(\tau)-G(\tau)N^{-1}(\tau)\beta(\tau))_H d\tau+\eta w(s)\label{eq:7-4}
\end{gather}
 We see easily that $X(s)$ can be written as
\begin{equation}\label{eq:7-5}
	X(s)=X_{0}(s)+\xi(s) \quad \text{ with } \quad \dfrac{d}{ds}\xi(s)=\mathcal{F}(s)\xi(s)+G(s)\Vt(s),\quad \xi(0)=0.
\end{equation}
 Clearly $\xi(s)\in\cH_1$. Moreover the functional $J(V(\cdot))$
in (\ref{eq:7-2}) can be written as 
\begin{multline}
J(V(\cdot))=\int_{0}^{T}\left[\dfrac{1}{2}(\mathcal{M}(s)X_{0}(s),X_{0}(s))_H -\dfrac{1}{2}(N^{-1}(s)\beta(s),\beta(s))_U +(\alpha(s),X_{0}(s))_H \right]ds\\
+\dfrac{1}{2}(\mathcal{M}_{T}X_{0}(T),X_{0}(T))_H +(\alpha_T,X_{0}(T))_H +\Jt(\xi(\cdot))\label{eq:7-7}
\end{multline}
 with 
\begin{multline}
\Jt(\xi(\cdot)):=\int_{0}^{T}\left[\dfrac{1}{2}(\mathcal{M}(s)\xi(s),\xi(s))_H +\dfrac{1}{2}(N(s)v(s),v(s))_U +(\mathcal{M}(s)X_{0}(s)+\alpha(s),\xi(s))_H \right]ds\\
+\dfrac{1}{2}(\mathcal{M}_{T}\xi(T),\xi(T))_H +(\mathcal{M}_{T}X_{0}(T)+\alpha_T,\xi(T))_H \label{eq:7-8}
\end{multline}

\subsection{KERNEL SOLUTION OF THE CONTROL PROBLEM }

In view of the definition of the norm in $\cH_1$ in (\ref{eq:6-3}) we can write 
\begin{equation}
\Jt(\xi(\cdot))=\dfrac{1}{2}\|\xi(\cdot)\|_{\cH_1}^{2}+\int_{0}^{T}(\mathcal{M}(s)X_{0}(s)+\alpha(s),\xi(s))_H ds+(\mathcal{M}_{T}X_{0}(T)+\alpha_T,\xi(T))_H \label{eq:7-9}
\end{equation}
and in view of the reproducing property (\ref{eq:6-5}) we get 
\begin{equation}
\Jt(\xi(\cdot))=\dfrac{1}{2}\|\xi(\cdot)\|_{\cH_1}^{2}+\left(\int_{0}^{T}K_1(\cdot,s)(\mathcal{M}(s)X_{0}(s)+\alpha(s))ds+K_1(\cdot,T)(\mathcal{M}_{T}X_{0}(T)+\alpha_T),\xi(\cdot)\right)_{\cH_1}\label{eq:7-10}
\end{equation}
The optimal $\hat{\xi}(\cdot)$ is immediately given by the formula 
\begin{equation}
\hat{\xi}(s)=-\int_{0}^{T}K_1(s,t)(\mathcal{M}(t)X_{0}(t)+\alpha(t))dt-K_1(s,T)(\mathcal{M}_{T}X_{0}(T)+\alpha_T)\label{eq:7-11}
\end{equation}
From the formula of the kernel given in (\ref{eq:6-6}), we obtain 
\begin{multline}
\hat{\xi}(s)=-\int_{0}^{T}\left(\int_{0}^{\min(s,t)}\Phi_{\cF,\cP}(s,\tau)G(\tau)N^{-1}(\tau)G^{*}(\tau)\Phi_{\cF,\cP}^{*}(t,\tau)\bE[\mathcal{M}(t)X_{0}(t)+\alpha(t)|\cW_{X}^{\tau}]d\tau\right)dt\\
-\int_{0}^{s}\Phi_{\cF,\cP}(s,\tau)G(\tau)N^{-1}(\tau)G^{*}(\tau)\Phi_{\cF,\cP}^{*}(T,\tau)\bE[\mathcal{M}_{T}X_{0}(T)+\alpha_T|\cW_{X}^{\tau}]d\tau\label{eq:7-12}
\end{multline}
Interchanging the order of integration in the first term at the right of (\ref{eq:7-12}), we derive that  
\begin{equation}
\hat{\xi}(s)=-\int_{0}^{s}\Phi_{\cF,\cP}(s,\tau)G(\tau)N^{-1}(\tau)G^{*}(\tau) \underbrace{\bE\left[\int_{\tau}^{T}\Phi_{\cF,\cP}^{*}(t,\tau)(\mathcal{M}(t)X_{0}(t)+\alpha(t))dt+\Phi_{\cF,\cP}^{*}(T,\tau)(\mathcal{M}_{T}X_{0}(T)+\alpha_T)|\cW_{X}^{\tau}\right]}_{=:r(\tau)} \label{eq:7-13}
\end{equation}
We now show that $r(s)$ is solution of a SBDE on $H$. Using that $\Phi_{\cF,\cP}$ is a semigroup, hence $\Phi_{\cF,\cP}(t,s)=\Phi_{\cF,\cP}(t,0)\Phi_{\cF,\cP}(0,s)$, we obtain that
\[
r(s)=\Phi_{\cF,\cP}^{*}(0,s)\bE\left[\int_{s}^{T}\Phi_{\cF,\cP}^{*}(t,0)(\mathcal{M}(t)X_{0}(t)+\alpha(t))dt+\Phi_{\cF,\cP}^{*}(T,0)(\mathcal{M}_{T}X_{0}(T)+\alpha_T)|\cW_{X}^{s}\right]
\]
hence
\begin{equation}
r(s)=\Phi_{\cF,\cP}^{*}(0,s)\left\{ -\int_{0}^{s}\Phi_{\cF,\cP}^{*}(t,0)(\mathcal{M}(t)X_{0}(t)+\alpha(t))dt+\bE\left[\Pi_T|\cW_{X}^{s}\right]\right\}\label{eq:7-17}
\end{equation}
where
\[
\Pi_T=\int_{0}^{T}\Phi_{\cF,\cP}^{*}(t,0)(\mathcal{M}(t)X_{0}(t)+\alpha(t))dt+\Phi_{\cF,\cP}^{*}(T,0)(\mathcal{M}_{T}X_{0}(T)+\alpha_T)
\]
Notice that there is no conditional expectation in the first term of \eqref{eq:7-17} because $t\in[0,s]$ so the content is within the filtration $\cW_{X}^{s}$. From the representation of stochastic integrals, we have
\begin{equation}
\bE[\Pi_T|\cW_{X}^{s}]=\bE[\Pi_T|\sigma(X^0)]+\int_{0}^{s}L_{X}(\tau)dw(\tau)\label{eq:7-18}
\end{equation}
 where the process $L_{X}\in L_{\cW_{X}}^{2}(0,T;H)$. We can
then write 
\begin{equation}
r(s)=-\int_{0}^{s}\Phi_{\cF,\cP}^{*}(t,s)(\mathcal{M}(t)X_{0}(t)+\alpha(t))dt+\Phi_{\cF,\cP}^{*}(0,s)\left(\bE[\Sigma|\sigma(X^0)]+\int_{0}^{s}L_{X}(\tau)dw(\tau)\right)\label{eq:7-20}
\end{equation}
 This implies 
\begin{equation}
-dr(s)=[(\mathcal{F}^{*}(s)-\cP(s)G(s)N^{-1}(s)G^{*}(s))r(s)+\mathcal{M}(s)X_{0}(s)+\alpha(s)]ds-\Gamma_{X}(s)dw(s), \quad r(T)=\mathcal{M}_{T}X_{0}(T)+\alpha_T\label{eq:7-21}
\end{equation}
with $\Gamma_{X}(\cdot)=\Phi_{\cF,\cP}^{*}(0,\cdot)L_{X}(\cdot)\in L_{\cW_{X}}^{2}(0,T;H)$. We define 
\begin{equation}
\hat X(s)=X_{0}(s)+\hat{\xi}(s).\label{eq:7-22}
\end{equation}
Using (\ref{eq:7-4}) and (\ref{eq:7-13}), we obtain after rearrangements
\begin{multline}
d\hat X(s)=[(\mathcal{F}(s)-G(s)N^{-1}(s)G^{*}(s)\cP(s))\hat X(s)+f(s)-G(s)N^{-1}(s)\beta(s)]ds+\eta dw(s)\\
+G(s)N^{-1}(s)G^{*}(s)(\cP(s)X_{0}(s)-r(s))ds, \quad \hat X(0)=X^0.\label{eq:7-23}
\end{multline}
 Calling $\lambda(s)=-\cP(s)X_{0}(s)+r(s)$ we obtain easily that it is the solution of the following BSDE
 \begin{multline*}
 	-d\lambda(s)=(\mathcal{F}^{*}(s)-\cP(s)G(s)N^{-1}(s)G^{*}(s))\lambda(s)ds+(\cP(s)(f(s)-G(s)N^{-1}(s)\beta(s))
 	+\alpha(s))ds\\+(\cP(s)\eta-\Gamma_{X}(s))dw(s), \quad \lambda(T)=\alpha_T.
 \end{multline*}
This linear BSDE has a unique solution. Since the terminal condition and all the terms of the drift part are deterministic, we can look for a solution such that $\Gamma_{X}(s)=\cP(s)\eta$, making the stochastic term vanish, and such that
\begin{equation}
-\dfrac{d}{ds}\lambda(s)=(\mathcal{F}^{*}(s)-\cP(s)G(s)N^{-1}(s)G^{*}(s))\lambda(s)+\cP(s)(f(s)-G(s)N^{-1}(s)\beta(s))+\alpha(s), \quad \lambda(T)=\alpha_T.\label{eq:7-24}
\end{equation}
This is exactly the function obtained in (\ref{eq:5-14}) through the method of completion of square.

 \bigskip
 
 \noindent{\bf Acknowledgements}
 
 \noindent
 The first author was funded by the European Research Council (grant REAL 947908). The second author was supported by the National Science Foundation under grant NSF-DMS-2204795.
 
 \bibliographystyle{plainnat}
 \bibliography{biblioPCAF,Kernels_MFC} 
 
 \subsection*{APPENDIX: EXPLICIT FORM OF THE KERNEL}
 The kernel $K$ was expressed in \eqref{eq:3-3} using the operators $\Phi_{\cF,\cP}$ and $\cP$ which act on $H$, so on both $X_x$ and its mean $\Xb$. Here we decouple these dependencies to obtain a more explicit form. We take 
 \begin{equation}
 	(\cJ_0 X)_{x}=J_0X_{x}+\Jb_0\Xb, \quad (\cJ_T X)_{x}=M_T X_{x}+\Mb_T\Xb  \label{eq:3-300}
 \end{equation}
 with $J_0,\Jb_0\in\cL(\R^{n},\R^{n})$ symmetric positive definite and $J$ invertible. The operator $\cP(s)$ can be written as 
 \begin{equation}
 	(\cP(s)X)_{x}=P(s)X_{x}+\Gamma(s)\Xb\label{eq:3-301}
 \end{equation}
 with $P(s)$ and $\Gamma(s)$ solutions of (\ref{eq:2-13}),(\ref{eq:2-16}). We furthermore have that
 \begin{gather*}
 	((\cJ_0 +\cP(0))X)_{x}=(J_0+P(0))X_{x}+(\Jb_0+\Gamma(0))\Xb \label{eq:3-305}\\
 	((\cJ_0 +\cP(0))^{-1}X)_{x}=(J_0+P(0))^{-1}X_{x}+((J_0+\Jb_0+\Sigma(0))^{-1}-(J_0+P(0))^{-1})\Xb \label{eq:3-306}
 \end{gather*}
 We next express the semigroup on $H$, denoted $\Phi_{\cF,\cP}(s,t)$, defined by $X(s)=\Phi_{\cF,\cP}(s,t)X(t),\:s>t$ with 
 \begin{equation}
 	\dfrac{d}{ds}X(s)=(\mathcal{F}(s)-G(s)N^{-1}(s)G^{*}(s)\cP(s))X(s)\label{eq:3-307}
 \end{equation}
 Since 
 \begin{equation*}
 	((\mathcal{F}(s)-G(s)N^{-1}(s)G^{*}(s)\cP(s))X)_{x}= (F(s)-G(s)N^{-1}(s)G^{*}(s)P(s))X_{x}+(\Fb(s)-G(s)N^{-1}(s)G^{*}(s)\Gamma(s))\Xb(s),
 \end{equation*}
 equation (\ref{eq:3-307}) reads 
 \begin{gather}
 	\dfrac{d}{ds}X_{x}(s)=(F(s)-G(s)N^{-1}(s)G^{*}(s)P(s))X_{x}(s)+(\Fb(s)-G(s)N^{-1}(s)G^{*}(s)\Gamma(s))\Xb(s)\label{eq:3-309}\\
 	\dfrac{d}{ds}\Xb(s)=(F(s)+\Fb(s)-G(s)N^{-1}(s)G^{*}(s)\Sigma(s))\Xb(s)\label{eq:3-310}
 \end{gather}
 We denote by $\Psi_{F,\Sigma}(s,t)$ the semigroup on $\R^{n}$ associated with $F(s)+\Fb(s)-G(s)N^{-1}(s)G^{*}(s)\Sigma(s)$, such that 
 \begin{equation*}
 	\Xb(s)=\Psi_{F,\Sigma}(s,t)X(t),\;s>t
 \end{equation*}
 From (\ref{eq:3-309}), (\ref{eq:3-310}) we obtain also 
 \[
 \dfrac{d}{ds}(X_{x}(s)-\Xb(s))=(F(s)-G(s)N^{-1}(s)G^{*}(s)P(s))(X_{x}(s)-\Xb(s))
 \]
 Therefore, if we denote by $\Psi_{F,P}(s,t)$ the semigroup on $\R^{n}$
 associated with $F(s)-G(s)N^{-1}(s)G^{*}(s)P(s)$, we can write 
 \begin{equation*}
 	X_{x}(s)-\Xb(s)=\Psi_{F,P}(s,t)(X_{x}(t)-\Xb(t))
 \end{equation*}
 Collecting results, we can write 
 \begin{equation*}
 	X_{x}(s)=\Psi_{F,P}(s,t)X_{x}(t)+(\Psi_{F,\Sigma}(s,t)-\Psi_{F,P}(s,t))\Xb(t)
 \end{equation*}
 So we have proven the formula
 \begin{equation}
 	(\Phi_{\cF,\cP}(s,t)X)_{x}=\Psi_{F,P}(s,t)X_{x}+(\Psi_{F,\Sigma}(s,t)-\Psi_{F,P}(s,t))\Xb\label{eq:3-114}
 \end{equation}
Concerning $K$, we have the following
 \begin{prop}
 	\label{prop3-10} The kernel $K$ defined in (\ref{eq:3-3}) is given by 
 	\begin{multline}\label{eq:3-116}
 		(K_1(s,t)X)_{x}=\Psi_{F,P}(s,0)(J+P(0))^{-1}\Psi_{F,P}^{*}(t,0)X_{x}+\int_{0}^{\min(t,s)}\Psi_{F,P}(s,\tau)G(\tau)N^{-1}(\tau)G^{*}(\tau)\Psi_{F,P}^{*}(t,\tau)d\tau\:X_{x}\\
 		+\left(\Psi_{F,\Sigma}(s,0)(J+\Jb+\Sigma(0))^{-1}\Psi_{F,\Sigma}^{*}(t,0)-\Psi_{F,P}(s,0)(J+P(0))^{-1}\Psi_{F,P}^{*}(t,0)\right)\Xb \\
 		+\int_{0}^{\min(t,s)}\left(\Psi_{F,\Sigma}(s,\tau)G(\tau)N^{-1}(\tau)G^{*}(\tau)\Psi_{F,\Sigma}^{*}(t,\tau)-\Psi_{F,P}(s,\tau)G(\tau)N^{-1}(\tau)G^{*}(\tau)\Psi_{F,P}^{*}(t,\tau)\right)d\tau\:\Xb.
 	\end{multline}
 \end{prop}
 
 \begin{proof}
 	We apply the definition (\ref{eq:3-3}) and use formulas (\ref{eq:3-306}) and (\ref{eq:3-114}).
 \end{proof}
 Expression \eqref{eq:3-116} enables direct computation of the kernel $K$ through the techniques exposed in the Appendix of \cite{aubin2022Kalman} based on matrix exponentials.
\end{document}